\documentclass[reqno]{amsart}
\usepackage[a4paper]{geometry}

\usepackage[T2A]{fontenc}
\usepackage[utf8]{inputenc}
\usepackage[ukrainian,russian,english]{babel}
\usepackage{comment}

\usepackage{amssymb}
\usepackage[mathscr]{eucal}
\usepackage{tikz}
\usetikzlibrary{positioning,calc,decorations.markings}
\usetikzlibrary{calc}
\usepackage{mathtools}
\usepackage{graphicx}
\usepackage{enumerate}
\usepackage{accents}
\usepackage{dsfont}
\usepackage{color}
\usepackage[colorlinks=true]{hyperref}
\hypersetup{urlcolor=blue,citecolor=red,linkcolor=blue}
\usepackage[initials]{amsrefs}
\definecolor{darkgreen}{rgb}{0.13, 0.55, 0.13}
\numberwithin{equation}{section}
\theoremstyle{plain}
\newtheorem{theorem}{Theorem}[section]
\newtheorem{proposition}[theorem]{Proposition}

\theoremstyle{definition}
\newtheorem*{rh-pb*}{Basic RH problem}
\newtheorem*{rh-I*}{RH problem normalized as $I$ at $\infty$}
\newtheorem*{sol-rh-pb*}{Soliton RH problem}
\newtheorem*{data*}{Data of this RH problem associated with $\BS{u_0(x)}$}
\theoremstyle{remark}
\newtheorem{remark}[theorem]{Remark}
\newtheorem*{notations*}{Notations}
\providecommand{\BS}[1]{\boldsymbol{#1}}  
\providecommand{\D}[1]{\mathbb{#1}}
\newcommand{\dd}{\mathrm{d}}
\newcommand{\eul}{\mathrm{e}}
\newcommand{\ee}{\mathrm{e}}
\newcommand{\ii}{\mathrm{i}}
\newlength{\dhatheight}
\newcommand{\doublehat}[1]{%
    \settoheight{\dhatheight}{\ensuremath{\hat{#1}}}%
    \addtolength{\dhatheight}{-0.3ex}%
    \hat{\vphantom{\rule{1pt}{\dhatheight}}%
    \smash{\hat{#1}}}}

\renewcommand{\Im}{\operatorname{Im}}
\renewcommand{\Re}{\operatorname{Re}}
\newcommand{\ord}{\mathrm{O}}
\DeclareMathOperator{\Res}{Res}

\newcommand{\red}{\textcolor{red}}

\newif\ifshort
\shorttrue

\title{The periodic Short-Pulse equation by the Riemann-Hilbert problem approach}
\date{December 2025}

\author[Iryna Karpenko]{Iryna Karpenko}
\address{IK: Faculty of Mathematics\\ University of Vienna\\
Oskar-Morgenstern-Platz 1\\ 1090 Wien\\ Austria\\ and B. Verkin Institute for Low Temperature Physics and Engineering\\ 47, Nauky ave\\ 61103 Kharkiv\\ Ukraine}
\email{\href{mailto:iryna.karpenko@univie.ac.at}{iryna.karpenko@univie.ac.at}}

\author[Dmitry Shepelsky]{Dmitry Shepelsky}
\address{DS: B.Verkin Institute for Low Temperature Physics and Engineering, Kharkiv, Ukraine}
\email{\href{mailto:shepelsky@yahoo.com }{shepelsky@yahoo.com}}

\author[Gerald Teschl]{Gerald Teschl}
\address{GT: Faculty of Mathematics\\ University of Vienna\\
Oskar-Morgenstern-Platz 1\\ 1090 Wien\\ Austria}
\email{\href{mailto:gerald.teschl@univie.ac.at}{gerald.teschl@univie.ac.at}}

\date{}

\begin{document}

\begin{abstract} In this paper, we formulate a Riemann--Hilbert problem framework --- based on the Unified Transform Method (Fokas method) --- for the Short Pulse equation posed on a periodic domain. Under the stated spectral and reconstruction regularity assumptions, we show that a local solution of the periodic initial-value problem can be reconstructed from the solution of a suitably defined Riemann–Hilbert problem, whose jump data are explicitly characterized by spectral functions associated with the initial condition.

\end{abstract}

\maketitle

\section{Introduction}

We consider the Short--Pulse (SP) equation
\begin{equation}\label{SP_1}
u_{xt} = u + \frac{1}{6}(u^{3})_{xx}.
\end{equation}
Setting $w=u_x$ and $q=\sqrt{1+w^2}$,  it can be rewritten as
\begin{equation}\label{SP-w-evolution}
w_t = u q^2 + \frac{u^2}{2}w_x.
\end{equation}
Every smooth solution of the short-pulse equation satisfies the following ``conservation law''
\begin{subequations}\label{SP}
\begin{align}\label{SP-2}
&q_t = \frac{1}{2}(u^2q)_x,\\
\label{tq}
&q=\sqrt{1+w^2},\\
\label{tw}
&w=u_x.
\end{align}
\end{subequations}
It was proposed by Sch\"afer and Wayne \cites{SW04,CJWS05} as an alternative to the nonlinear Schr\"odinger (NLS) equation
for describing the propagation of ultra--short, intense optical pulses.
Unlike the NLS equation, which models slow amplitude modulation under narrowband assumptions,
the SP equation is applicable in regimes where these assumptions fail, as encountered in
high--speed fiber--optic communication.
Numerical simulations in \cite{CJWS05} confirmed that the SP equation provides a reliable
description of pulse dynamics in the presence of broad spectral bandwidths.

The SP equation is also closely related to the Camassa--Holm (CH) equation and can
be viewed as a short--wave limit of the modified Camassa--Holm equation, also known
as the Fuchssteiner--Olver--Rosenau--Qiao equation
\cites{F96,FO81,Q06,OR96}. This structural connection motivates the application of
analytical techniques developed for CH-type equations to the study of the SP
equation.

The SP equation has been studied extensively. Solitary-wave solutions
were obtained in \cites{Gambino2015,Matsuno2007,Sakovich2006}; periodic
solutions were constructed in \cite{Matsuno2008}. The well-posedness of the Cauchy problem was
investigated in \cites{Coclite2015,Pelinovsky2010,SW04}, while sufficient
conditions for wave breaking were derived in \cite{Liu2009}.

We consider the SP equation on a finite spatial interval \( x \in [0, L] \), supplemented with the periodicity conditions
\begin{equation}\label{period}
    u(0,t) = u(L,t), \quad u_x(0,t) = u_x(L,t),
\end{equation}
and the initial condition
\begin{equation}\label{ic}
    u(x,0) = u_0(x), \quad x \in [0,L].
\end{equation}
We assume that the initial data are real-valued, sufficiently smooth, and periodic, and satisfy
\begin{equation}\label{mean-zero-data}
    \int_0^L u_0(x)\dd x=0.
\end{equation}
Indeed, integrating \eqref{SP_1} over one period shows that every smooth periodic solution has zero mean. To fix ideas, we additionally assume that $u_{0x}(0)=u_{0x}(L)=0$ (cf. \cite{KS25_}). 

In this context, notice that, due to periodicity, we can shift the spatial variable. Since $u_0$ is differentiable and periodic, Rolle's theorem gives a point $x_0$ such that $u_x(x_0,0)=0$. Thus, we can consider $\tilde u(x,t)=u(x+x_0,t)$. On the one hand, $\tilde u$ satisfies the same equation as $u$. On the other hand, $\tilde u_x(0,0)=u_x(x_0,0)=0$ and $\tilde u_x(L,0)=u_x(x_0+L,0)=0$.

To investigate this problem, we use the inverse scattering transform (IST). The IST is a powerful analytical method to solve a class of one-dimensional nonlinear evolution equations known as integrable systems. Introduced more than fifty years ago, it is based on the analysis of an associated Lax pair whose compatibility condition yields the underlying nonlinear equation \cites{IST01,IST02,ablowitz1974inverse,ablowitz1991solitons}. The approach involves three main steps: the solution of a direct spectral problem mapping initial data to spectral data, the determination of the time evolution of these spectral data, and the reconstruction of the solution through an inverse problem. In this sense, the IST may be viewed as a nonlinear analog of the Fourier transform.

For initial value problems on the real line with rapidly decaying data, the inverse step can be carried out through the Gelfand–Levitan–Marchenko integral equations \cite{novikov1984theory}, or equivalently formulated as a Riemann–Hilbert (RH) factorization problem in the complex spectral plane \cite{trogdon2015riemann}. Extending the IST beyond this setting to problems posed on the half-line or on finite intervals presents additional analytical challenges due to the presence of boundary conditions.

These difficulties were overcome by the Unified Transform Method, introduced by Fokas \cites{fokas1997,fokas2008unified}, which analyzes both equations of the Lax pair simultaneously. A key feature of the method is the \emph{global relation}, which links the spectral transforms of the initial and boundary data. For certain boundary conditions, referred to as \emph{linearizable}, the global relation can be resolved explicitly, leading to a well-posed RH problem formulated solely in terms of prescribed data.

The periodic initial-boundary value problem for the nonlinear Schrödinger (NLS)
equation was shown to be linearizable in \cites{DFL21,FL21}, and finite-band
(algebro-geometric) solutions were characterized, based on this method, in \cite{SKPB24}. More recently, an analogous linearizability result was
established for the periodic IBV problem for the Camassa--Holm (CH) equation in \cite{BKSZ25}.

In the present paper, we develop the RH method for solving the initial value problem 
for the SP equation with periodic boundary
conditions. From a physical perspective, periodic boundary conditions may be viewed as an idealized model for pulse propagation in closed or spatially repeating optical structures, such as fiber loops or resonant systems, where a pulse undergoes repeated propagation through the same medium.

Correspondingly, our main aim 
is to obtain a representation of the solution 
of this problem in terms of the solution of an associated RH problem, the data for which
are determined in terms of the spectral functions associated with the initial data.

The outline of the paper is as follows.
In Section \ref{sec:2} we introduce the eigenfunctions and spectral functions near $k=\infty$ and $k=0$. Section \ref{sec:3} presents the global relations, and Section \ref{sec:4} formulates the pre-Riemann--Hilbert problem on an interval. In Section \ref{sec:5} we introduce the $(y,t)$ variables and describe the reconstruction of a local solution. Section \ref{sec:6} gives the master RH problem for the periodic problem and verifies the initial conditions and periodicity. 

\begin{notations*}
    Throughout, $f^*(k):=\overline{f(\bar k)}$ for a scalar spectral function $f$, $\mathbb C^\pm:=\{k:\pm\Im k>0\}$, and
\[
 \sigma_2=\begin{pmatrix}0&-\ii\\\ii&0\end{pmatrix},
 \qquad
 \sigma_3=\begin{pmatrix}1&0\\0&-1\end{pmatrix}.
\]
We use $\eul^{z\hat\sigma_3}A:=\eul^{z\sigma_3}A\eul^{-z\sigma_3}$; a dot over a spectral function denotes differentiation with respect to $k$.
\end{notations*}

\section{Eigenfunctions and Spectral Functions}\label{sec:2}

The Short Pulse equation is formally integrable: it is the compatibility condition for the Lax pair (\cite{SS05}):
\begin{subequations}\label{Lax}
\begin{align}\label{Lax-x}
\Phi_x&=U\Phi,\\
\label{Lax-t}
\Phi_t&=V\Phi,
\end{align}
\end{subequations}
where the coefficients $U\equiv U(x,t,k)$ and $V\equiv V(x,t,k)$ are defined by
\begin{subequations}\label{Lax-UV}
\begin{align}
    U&=-\ii k\begin{pmatrix}
        1& u_x\\
         u_x&-1 
    \end{pmatrix},\\
        V&=\begin{pmatrix}
        -\frac{\ii k}{2}u^2+\frac{\ii}{4  k}&-\frac{\ii k}{2}u^2u_x-\frac{1}{2}u\\
        -\frac{\ii k}{2}u^2u_x+\frac{1}{2}u&\frac{\ii k}{2}u^2-\frac{\ii}{4k} 
    \end{pmatrix}.
\end{align} 
\end{subequations}

\subsection{Eigenfunctions near $k=\infty$}

In order to control the large $k$ behavior, we introduce 
\[
\tilde \Phi_\infty(x,t,k):= P(x,t)\Phi(x,t,k),
\]
where
\[
P(x,t):=\sqrt{\frac{1+q(x,t)}{2q(x,t)}}\begin{pmatrix}
1 & \frac{w(x,t)}{1+q(x,t)}  \\
-\frac{w(x,t)}{1+q(x,t)} & 1 \\
\end{pmatrix},
\]
where $w(x,t)=u_x(x,t)$ and $q(x,t)=\sqrt{1+w^2(x,t)}$. This transforms \eqref{Lax} into

\begin{subequations}\label{Lax-Q-form}
   \begin{align}
    \tilde \Phi_{\infty x}+\ii k p_x \sigma_3 \tilde \Phi_{\infty}&=U_\infty \tilde \Phi_{\infty },\\
    \tilde \Phi_{\infty t}+\ii k p_t \sigma_3 \tilde \Phi_{\infty}&=V_\infty \tilde \Phi_{\infty }
\end{align} 
\end{subequations}
where $U_\infty\equiv U_\infty(x,t,k)$ and $V_\infty\equiv V_\infty(x,t,k)$ are given by 
\begin{equation}\label{Uinf}
    U_\infty(x,t,k)=\frac{u_{xx}}{2q^2}\begin{pmatrix}
        0 & 1 \\
-1 & 0
    \end{pmatrix},
\end{equation}
and
\begin{equation}\label{Vinf}
    V_\infty(x,t,k)= -\frac{1}{4\ii k q}\begin{pmatrix}
        1-q & -w \\
-w & -1+q
    \end{pmatrix}+\frac{u^2u_{xx}}{4q^2}\begin{pmatrix}
        0 & 1 \\
-1 & 0
    \end{pmatrix}
\end{equation}
and $p(x,t,k)$ is defined by 
\begin{equation}\label{p}
p(x,t,k)=\int_0^x q(\xi,t)\dd\xi+\int_0^t\frac{u^2(0,\zeta)q(0,\zeta)}{2} \dd\zeta-\frac{t}{4k^2}.
\end{equation}

\begin{remark}
     The SP equation \eqref{SP} implies the following ``conservation law'' in the periodic case:
\begin{equation}\label{cons_law}
    \int_0^L q(\xi,t)\dd\xi= \int_0^L q(\xi,0)\dd\xi:=\theta.
\end{equation}
     
\end{remark}

Notice that $p(x,t,k)$ satisfies
\[
p_x=q,\qquad
p_t=\frac{u^2q}{2}-\frac{1}{4k^2},\qquad p(0,0,k)=0,\qquad p(L,0,k)=\theta.
\]

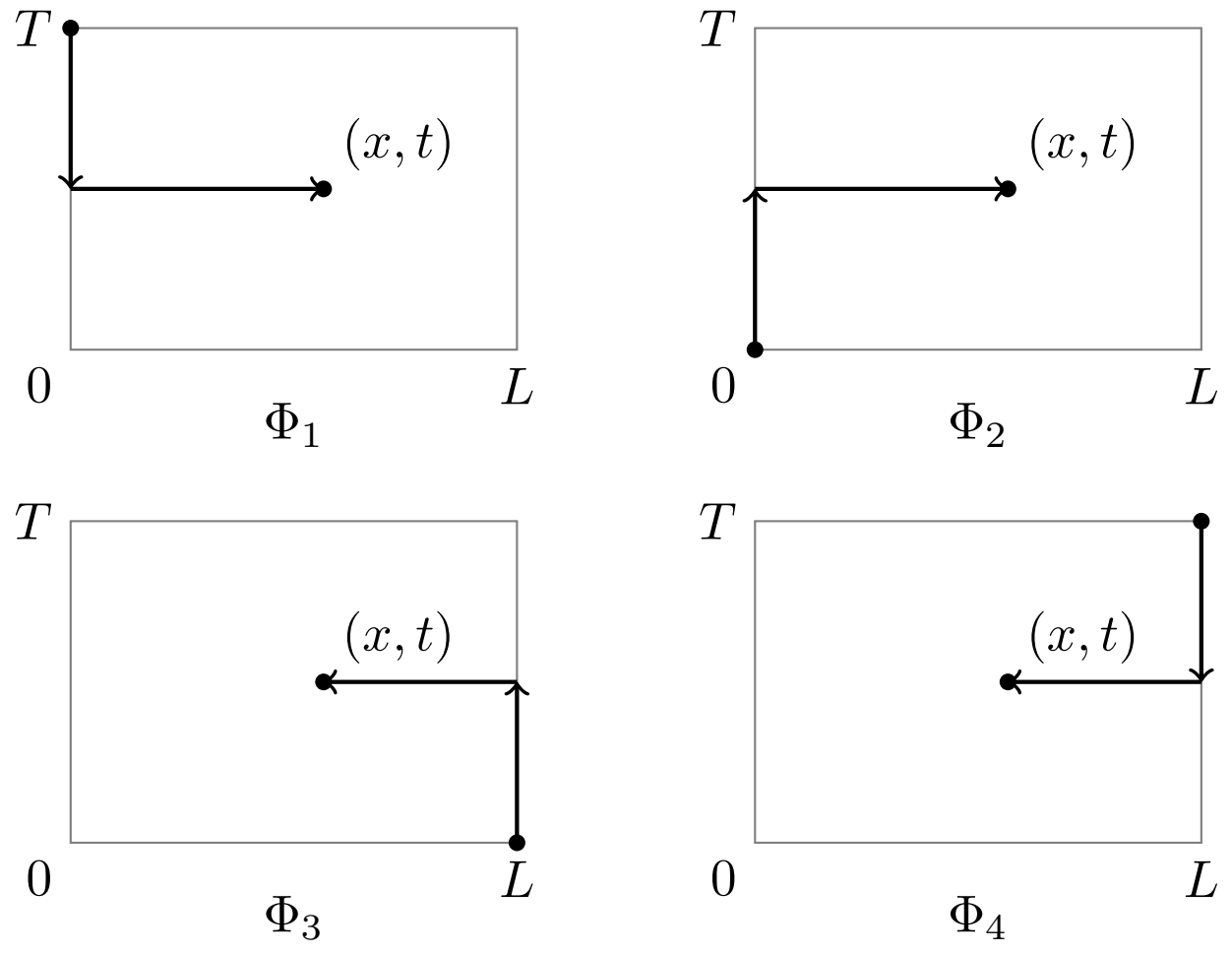
\begin{figure}[h]
    \centering
    \begin{tikzpicture}[x=1cm,y=.72cm]
 \foreach \j/\dx/\dy/\sx/\sy in
 {1/0/4.6/0/3,2/4.6/4.6/0/0,3/0/0/3/0,4/4.6/0/3/3}
 {
  \begin{scope}[shift={(\dx,\dy)}]
   \draw[gray] (0,0) rectangle (3,3);
   \node[below left] at (0,0) {$0$};
   \node[below] at (3,0) {$L$};
   \node[left] at (0,3) {$T$};
   \draw[thick,->] (\sx,\sy)--(\sx,1.5);
   \draw[thick,->] (\sx,1.5)--(1.7,1.5);
   \fill (\sx,\sy) circle (1.6pt);
   \fill (1.7,1.5) circle (1.6pt);
   \node[above right] at (1.7,1.5) {$(x,t)$};
   \node at (1.5,-.7) {$\Phi_{\j}$};
  \end{scope}
 }
\end{tikzpicture}
    \caption{Paths of integration for $\Phi_j$}
    \label{fig:paths}
\end{figure}

Introduce the solutions $\tilde\Phi_{\infty j}$, $j=1,2,3,4$ to \eqref{Lax-Q-form} through the solutions $\Phi_{\infty j}$ of the integral equations, which are fixed by the associated
initial 
integration point. Namely, $\tilde \Phi_{\infty }(x,t,k)=\Phi_{\infty }(x,t,k)\eul^{-\ii k p(x,t,k)\sigma_3}$, where $\Phi_{\infty j}$, $j=1,2,3,4$ satisfy
the integral equation 
\begin{equation}\label{inteq_inf}
\Phi_{\infty}(x,t,k)=I+\int_{(x^*,t^*)}^{(x,t)}
	\eul^{-ik(p(x,t,k)-p(y,\tau,k))\hat\sigma_3}( U_\infty\Phi_\infty \dd y+V_\infty\Phi_\infty \dd \tau)(y,\tau,k),
\end{equation}
where $(x^*,t^*)$ is taken respectively as $(0,T)$, $(0,0)$, $(L,0)$, $(L,T)$ (with some $T>0$). The integration paths are those shown in Figure~\ref{fig:paths}; due to the compatibility of the Lax pair, the value of the integral is independent of the chosen path.

 More precisely,
\begin{equation}\label{inteq_inf1}
    \begin{aligned}
\Phi_{\infty 1}(x,t,k)=&I+
\int_{0}^{x}
	\eul^{-ik\int_y^xq(\xi,t)\dd \xi\hat\sigma_3}( U_\infty\Phi_{\infty1})(y,t,k) \dd y  - \\
&\eul^{-ik\int_0^xq(\xi,t)\dd \xi\hat\sigma_3}
\int_{t}^{T}
	\eul^{-ik\left(\int_\tau^t \frac{u^2(0,\zeta)q(0,\zeta)}{2}\dd \zeta-\frac{1}{4k^2}(t-\tau)\right)\hat\sigma_3}( V_\infty\Phi_{\infty1})(0,\tau,k) \dd \tau,
\end{aligned}
\end{equation}
    
\begin{equation}\label{inteq_inf2}
     \begin{aligned}
\Phi_{\infty 2}(x,t,k)=&I+
\int_{0}^{x}
	\eul^{-ik\int_y^xq(\xi,t)\dd \xi\hat\sigma_3}( U_\infty\Phi_{\infty2})(y,t,k) \dd y  + \\
&\eul^{-ik\int_0^xq(\xi,t)\dd \xi\hat\sigma_3}
\int_{0}^{t}
	\eul^{-\ii k\left(\int_\tau^t \frac{u^2(0,\zeta)q(0,\zeta)}{2}\dd \zeta-\frac{1}{4k^2}(t-\tau)\right)\hat\sigma_3}( V_\infty\Phi_{\infty2})(0,\tau,k) \dd \tau,
\end{aligned}
\end{equation}

\begin{equation}\label{inteq_inf3}
        \begin{aligned}
\Phi_{\infty 3}(x,t,k)=&I-
\int_{x}^{L}
	\eul^{-ik\int_y^x q(\xi,t)\dd \xi\hat\sigma_3}( U_\infty\Phi_{\infty3})(y,t,k) \dd y  + \\
&\eul^{ik\int_x^Lq(\xi,t)\dd \xi\hat\sigma_3}
\int_{0}^{t}
	\eul^{-\ii k\left(\int_\tau^t \frac{u^2(0,\zeta)q(0,\zeta)}{2}\dd \zeta-\frac{1}{4k^2}(t-\tau)\right)\hat\sigma_3}( V_\infty\Phi_{\infty3})(L,\tau,k) \dd \tau,
\end{aligned}
\end{equation}
   
\begin{equation}\label{inteq_inf4}
        \begin{aligned}
\Phi_{\infty 4}(x,t,k)=&I-
\int_{x}^{L}
	\eul^{-ik\int_y^x q(\xi,t)\dd \xi\hat\sigma_3}( U_\infty\Phi_{\infty4})(y,t,k) \dd y  - \\
&\eul^{ik\int_x^Lq(\xi,t)\dd \xi\hat\sigma_3}
\int_{t}^{T}
	\eul^{-\ii k\left(\int_\tau^t \frac{u^2(0,\zeta)q(0,\zeta)}{2}\dd \zeta-\frac{1}{4k^2}(t-\tau)\right)\hat\sigma_3}( V_\infty\Phi_{\infty4})(L,\tau,k) \dd \tau,
\end{aligned}
\end{equation}

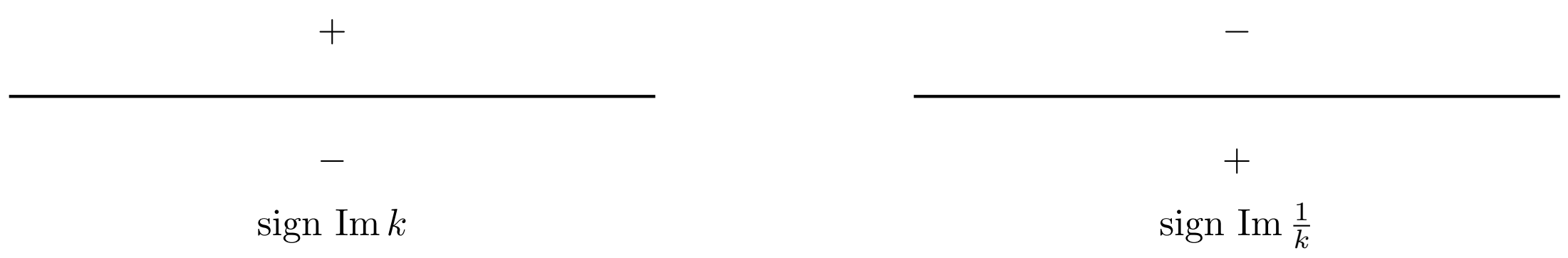
\begin{figure}[ht]
    \centering
\begin{tikzpicture}[scale=1.2]
        \node at (0,-1) {$\mathrm{sign}\,\Im k$};
    \draw[thick] (-2.5,0) -- (2.5,0);
    \node at (0,0.5) {$+$};
    \node at (0,-0.5) {$-$};
    \node at (7,-1) {$\mathrm{sign}\,\Im\frac{1}{k}$};
    \draw[thick] (4.5,0) -- (9.5,0);
    \node at (7,0.5) {$-$};
    \node at (7,-0.5) {$+$};
\end{tikzpicture}
    \caption{Signature tables for $p(x,t,k)$  }
    \label{fig:sign}
\end{figure}

The domains where the exponentials are bounded in the complex $k$-plane are separated by the real axis, see 
 Figure \ref{fig:sign}.

Let $A^{(1)}$ and $A^{(2)}$ denote the columns of a $2\times 2$ matrix $A = \bigl( A^{(1)}\ \ A^{(2)} \bigr)$. With this notation, using Neumann series expansions, we obtain the following properties of $\Phi_{\infty i}^{(j)}(x,t,k)$:

\begin{enumerate}
    \item $\Phi_{\infty 1}^{(1)}(x,t,k)$ is analytic in $\mathbb{C}\setminus \{0\}$.
 Moreover, $\Phi_{\infty 1}^{(1)}(0,t,k)$=$\begin{pmatrix}
        1\\0
    \end{pmatrix}+O(\frac{1}{k})$ as $k\to\infty$ in $\mathbb{C}^-$; $\Phi_{\infty 1}^{(1)}(x,T,k)$=$\begin{pmatrix}
        1\\0
    \end{pmatrix}+O(\frac{1}{k})$ as $k\to\infty$ in $\mathbb{C}^+$.

    \item $\Phi_{\infty 1}^{(2)}(x,t,k)$ is analytic in $\mathbb{C}\setminus \{0\}$.
 Moreover, $\Phi_{\infty 1}^{(2)}(0,t,k)$=$\begin{pmatrix}
        0\\1
    \end{pmatrix}+O(\frac{1}{k})$ as $k\to\infty$ in $\mathbb{C}^+$; $\Phi_{\infty 1}^{(2)}(x,T,k)$=$\begin{pmatrix}
        0\\1
    \end{pmatrix}+O(\frac{1}{k})$ as $k\to\infty$ in $\mathbb{C}^-$.

     \item $\Phi_{\infty 2}^{(1)}(x,t,k)$ is analytic in $\mathbb{C}\setminus\{0\}$ and $\Phi_{\infty 2}^{(1)}(x,0,k)$ is analytic in $\mathbb{C}$.
 Moreover, $\Phi_{\infty 2}^{(1)}(x,t,k)$=$\begin{pmatrix}
        1\\0
    \end{pmatrix}+O(\frac{1}{k})$ as $k\to\infty$ in $\mathbb{C}^+$.

    \item $\Phi_{\infty 2}^{(2)}(x,t,k)$ is analytic in $\mathbb{C}\setminus\{0\}$ and $\Phi_{\infty 2}^{(2)}(x,0,k)$ is analytic in $\mathbb{C}$.
 Moreover, $\Phi_{\infty 2}^{(2)}(x,t,k)$=$\begin{pmatrix}
        0\\1
    \end{pmatrix}+O(\frac{1}{k})$ as $k\to\infty$ in $\mathbb{C}^-$.

 \item $\Phi_{\infty 3}^{(1)}(x,t,k)$ is analytic in $\mathbb{C}\setminus \{0\}$ and $\Phi_{\infty 3}^{(1)}(x,0,k)$ is analytic in $\mathbb{C}$.
 Moreover, $\Phi_{\infty 3}^{(1)}(L,t,k)$=$\begin{pmatrix}
        1\\0
    \end{pmatrix}+O(\frac{1}{k})$ as $k\to\infty$ in $\mathbb{C}^+$; $\Phi_{\infty 3}^{(1)}(x,0,k)$=$\begin{pmatrix}
        1\\0
    \end{pmatrix}+O(\frac{1}{k})$ as $k\to\infty$ in $\mathbb{C}^-$.

    \item $\Phi_{\infty 3}^{(2)}(x,t,k)$ is analytic in $\mathbb{C}\setminus \{0\}$ and $\Phi_{\infty 3}^{(2)}(x,0,k)$ is analytic in $\mathbb{C}$.
 Moreover, $\Phi_{\infty 3}^{(2)}(L,t,k)$=$\begin{pmatrix}
        0\\1
    \end{pmatrix}+O(\frac{1}{k})$ as $k\to\infty$ in $\mathbb{C}^-$; $\Phi_{\infty 3}^{(2)}(x,0,k)$=$\begin{pmatrix}
        0\\1
    \end{pmatrix}+O(\frac{1}{k})$ as $k\to\infty$ in $\mathbb{C}^+$.

     \item $\Phi_{\infty 4}^{(1)}(x,t,k)$ is analytic in $\mathbb{C}\setminus \{0\}$.
 Moreover, $\Phi_{\infty 4}^{(1)}(x,t,k)$=$\begin{pmatrix}
        1\\0
    \end{pmatrix}+O(\frac{1}{k})$ as $k\to\infty$ in $\mathbb{C}^-$.

    \item $\Phi_{\infty 4}^{(2)}(x,t,k)$ is analytic in $\mathbb{C}\setminus \{0\}$.
 Moreover, $\Phi_{\infty 4}^{(2)}(x,t,k)$=$\begin{pmatrix}
        0\\1
    \end{pmatrix}+O(\frac{1}{k})$ as $k\to\infty$ in $\mathbb{C}^+$.

\end{enumerate}

Since the coefficient matrices in the Lax pair \eqref{Lax-Q-form} are traceless, we first note that
\[
\det \Phi_{\infty j} \equiv 1, \qquad j=1,2,3,4.
\]

Moreover, the matrices $U_\infty(k)\equiv U_\infty(x,t,k)$ and 
$V_\infty(k)\equiv V_\infty(x,t,k)$ satisfy the symmetries
\begin{subequations}\label{sym-UV}
\begin{alignat}{3}\label{sym-U}
U_\infty(- k)&=\overline{U_\infty(\bar k)},&\qquad&
U_\infty(-k)=\sigma_2U_\infty(k)\sigma_2,\\
V_\infty(- k)&=\overline{V_\infty(\bar k)},&&
V_\infty(-k)=\sigma_2V_\infty(k)\sigma_2
\end{alignat}
\end{subequations}
and the phase function $p(k)\equiv p(x,t,k)$ obeys
\begin{equation}\label{sym-p}
p^*(k)=p(k)=p(-k).
\end{equation}
These relations imply that the functions $\Phi_{\infty j}$ inherit the same symmetries:
\begin{equation}\label{sym-Phi}
\Phi_{\infty j}(k)=\overline{\Phi_{\infty j}(-\bar k)},\qquad
\Phi_{\infty j}(k)=\sigma_2\overline{\Phi_{\infty j}(\bar k)}\sigma_2.
\end{equation}

The functions $\tilde\Phi_{\infty j}$,
being the solutions of both equations of \eqref{Lax-Q-form}, are related (where they are defined) by matrices independent of $x$ and $t$, and thus
\begin{subequations}\label{rel_inf}
    \begin{align}\label{rel_inf_1}
      &\Phi_{\infty 3}(x,t,k)=\Phi_{\infty 2}(x,t,k)\eul^{-ikp(x,t,k)\sigma_3}s(k)\eul^{ikp(x,t,k)\sigma_3},\\\label{rel_inf_2}
      &\Phi_{\infty 1}(x,t,k)=\Phi_{\infty 2}(x,t,k)\eul^{-ikp(x,t,k)\sigma_3}S(k)\eul^{ikp(x,t,k)\sigma_3},\\\label{rel_inf_3}
      &\Phi_{\infty 4}(x,t,k)=\Phi_{\infty 3}(x,t,k)\eul^{-ikp(x,t,k)\sigma_3}S_1(k)\eul^{ikp(x,t,k)\sigma_3}.
    \end{align}
\end{subequations}
Due to the symmetries of $U_\infty$ and $V_\infty$,  $s(k)$ and $S(k)$
can be written as
\begin{subequations}
    \begin{align}\label{s}
       & s(k)=\Phi_{\infty 3}(0,0,k)=\begin{pmatrix}
       a^*( k)& b(k)\\
        -b^*( k) & a(k)
    \end{pmatrix},\\\label{S}
    & S(k)=\Phi_{\infty 1}(0,0,k)=\begin{pmatrix}
        A^*( k)& B(k)\\
        -B^*( k) & A(k)
    \end{pmatrix},\\\label{S_1}
& S_1(k)=\eul^{ik\theta\sigma_3}\Phi_{\infty 4}(L,0,k)\eul^{-ik\theta\sigma_3}=\begin{pmatrix}
        A_1^*( k)& \eul^{2ik\theta}B_1(k)\\
        -\eul^{-2ik\theta}B_1^*( k) & A_1(k)
    \end{pmatrix}.
    \end{align}
\end{subequations}

Considering  equation \eqref{inteq_inf3} for $t=0$ allows to derive the following properties of $a(k)$ and $b(k)$:

\begin{enumerate}
    \item $a(k)$ and $ b(k)$ are analytic in $\mathbb{C}\setminus\{0\}$. Moreover,  as $k\to\infty$ 
    \begin{subequations}\label{s_inf}
\begin{align}
a(k)&=1+O(\frac{1}{k})+O(\frac{\eul^{2\ii k\theta}}{k}),\\
a^*(k)&=1+O(\frac{1}{k})+O(\frac{\eul^{-2\ii k\theta}}{k}),\\
b(k)&=O(\frac{1}{k})+O(\frac{\eul^{2\ii k\theta}}{k}),\\
b^*(k)&=O(\frac{1}{k})+O(\frac{\eul^{-2\ii k\theta}}{k})
\end{align}
 \end{subequations}

\item The symmetries:
\begin{equation}\label{symm_a}
    a(-k)=\overline{a(\bar k)},~ b(-k)=\overline{b(\bar k)}
\end{equation}

\item Determinant relation:

\begin{equation}\label{det_rel}
    a(k)a^*(k)+ b(k)b^*(k)=1
\end{equation}
 
\end{enumerate}

Considering  equation \eqref{inteq_inf1} for $x=0$, we obtain the following properties of $A(k)$ and $B(k)$:

\begin{enumerate}

 \item $ A(k)$ and $ B(k)$ are analytic in $\mathbb{C}\setminus\{0\}$. Moreover, as $k\to\infty$
 \begin{subequations}
\label{S__inf}
\begin{align}
A(k)&=1+O(\frac{1}{k}),\quad k\in \mathbb{C}^+,\\
A^*(k)&=1+O(\frac{1}{k}),\quad k\in \mathbb{C}^-,\\
B(k)&=O(\frac{1}{k}),\quad k\in \mathbb{C}^+,\\
B^*(k)&=O(\frac{1}{k}),\quad k\in \mathbb{C}^-
\end{align}
 \end{subequations}
 
    \item The symmetries:
    \begin{equation}\label{symm_A}
    A(-k)=\overline{A(\bar k)},~ B(-k)=\overline{B(\bar k)}
\end{equation}

\item Determinant relation:
\begin{equation}\label{det_rel_A}
    A(k)\overline{A(\bar k)}+ B(k)\overline{B(\bar k)}=1
\end{equation}

\end{enumerate}

Similarly, considering  equation \eqref{inteq_inf4} for $x=L$, we obtain the following properties of $A_1(k)$ and $B_1(k)$:
\begin{enumerate}
    \item $ A_1(k)$ and $ B_1(k)$ are analytic in $\mathbb{C}\setminus\{0\}$. Moreover, as $k\to\infty$
    \begin{subequations}
        \label{S_1_inf}
\begin{align}
A_1(k)&=1+O(\frac{1}{k}),\quad k\in \mathbb{C}^+,\\
A_1^*(k)&=1+O(\frac{1}{k}),\quad k\in \mathbb{C}^-,\\
B_1(k)&=O(\frac{1}{k}),\quad k\in \mathbb{C}^+,\\
B_1^*(k)&=O(\frac{1}{k}),\quad k\in \mathbb{C}^-
\end{align}
    \end{subequations}

\item The symmetries:
\begin{equation}\label{symm_A1}
    A_1(-k)=\overline{A_1(\bar k)},~ B_1(-k)=\overline{B_1(\bar k)}
\end{equation}

\item Determinant relation:

\begin{equation}\label{det_rel_A1}
    A_1(k)\overline{A_1(\bar k)}+ B_1(k)\overline{B_1(\bar k)}=1
\end{equation}

\end{enumerate}

\subsection{Eigenfunctions near $k=0$}

In order to have good control as $k \to 0$, 
it is convenient to rewrite the original Lax pair \eqref{Lax} in the form

\begin{align}\label{Lax_0}
    \tilde \Phi_{0 x}+\ii k  \sigma_3 \tilde \Phi_{0}&=U_0 \tilde \Phi_{0 },\\
    \tilde \Phi_{0 t}+ \frac{1}{4\ii k}  \sigma_3 \tilde \Phi_{0}&=V_0 \tilde \Phi_{0 }
\end{align}
where
\begin{equation}\label{U0}
    U_0(x,t,k)=-\ii kw\begin{pmatrix}
        0 & 1 \\
1 & 0
    \end{pmatrix}
\end{equation}
and
\begin{align}\label{V0}
    V_0(x,t,k)=-&\frac{\ii k u^2}{2}\begin{pmatrix}
        1 & w \\
w& -1
    \end{pmatrix}+\frac{u}{2}\begin{pmatrix}
        0 & -1 \\
1 & 0
    \end{pmatrix}.
\end{align}
Notice that $U_0(x,t,0)\equiv 0$. 

Introduce (compare with \eqref{p})
\begin{equation}\label{p0}
p_0(x,t,k)=x-\frac{t}{4k^2},
\end{equation}
and 
the solutions $\tilde{\Phi}_{0j}(x,t,k)$, $j=1,2,3,4$, of 
\eqref{Lax_0} in the same way as $\tilde{\Phi}_{\infty j}(x,t,k)$ by defining
\[
\Phi_{0}(x,t,k):=\tilde \Phi_{0}(x,t,k)\eul^{\ii k p_0(x,t,k)\sigma_3},
\]
where the functions $\Phi_{0j}(x,t,k)$ are defined as the solutions of the 
integral equations
\begin{equation}\label{inteq_0}
\Phi_{0}(x,t,k)=I+\int_{(x^*,t^*)}^{(x,t)}
	\eul^{-ik\left((x-y)-\frac{t-\tau}{4 k^2}\right)\hat\sigma_3}( U_0\Phi_0 \dd y+V_0\Phi_0 \dd \tau)(y,\tau,k),
\end{equation}
Namely,
\begin{equation}\label{inteq_01}
        \begin{aligned}
\Phi_{0 1}(x,t,k)=&I+
\int_{0}^{x}
	\eul^{-ik(x-y)\hat\sigma_3}( U_0\Phi_{01})(y,t,k) \dd y  - \\
&\eul^{-ikx\hat\sigma_3}
\int_{t}^{T}
	\eul^{(\frac{\ii }{4k}(t-\tau))\hat\sigma_3}( V_0\Phi_{01})(0,\tau,k) \dd \tau,
\end{aligned}
\end{equation}

    \begin{equation}\label{inteq_02}
        \begin{aligned}
\Phi_{0 2}(x,t,k)=&I+
\int_{0}^{x}
	\eul^{-ik(x-y)\hat\sigma_3}( U_0\Phi_{02})(y,t,k) \dd y  + \\
&\eul^{-ikx\hat\sigma_3}
\int_{0}^{t}
	\eul^{(\frac{\ii }{4k}(t-\tau))\hat\sigma_3}( V_0\Phi_{02})(0,\tau,k) \dd \tau,
\end{aligned}
\end{equation}

    \begin{equation}\label{inteq_03}
        \begin{aligned}
\Phi_{0 3}(x,t,k)=&I-
\int_{x}^{L}
	\eul^{-ik(x-y)\hat\sigma_3}( U_0\Phi_{03})(y,t,k) \dd y  + \\
&\eul^{-ik(x-L)\hat\sigma_3}
\int_{0}^{t}
	\eul^{(\frac{\ii }{4k}(t-\tau))\hat\sigma_3}( V_0\Phi_{03})(L,\tau,k) \dd \tau,
\end{aligned}
\end{equation}

    \begin{equation}\label{inteq_04}
      \begin{aligned}
\Phi_{04}(x,t,k)=&I-
\int_{x}^{L}
	\eul^{-ik(x-y)\hat\sigma_3}( U_0\Phi_{04})(y,t,k) \dd y  - \\
&\eul^{-ik(x-L)\hat\sigma_3}
\int_{t}^{T}
	\eul^{(\frac{\ii }{4k}(t-\tau))\hat\sigma_3}( V_0\Phi_{04})(L,\tau,k) \dd \tau,
\end{aligned}
\end{equation}

Notice that $U_0(x,t,0)=0$.

Using Neumann series expansions, we conclude that $\Phi_{0 i}^{(j)}(x,t,k)$ are analytic in $\mathbb{C}\setminus \{0\}$. Moreover,
\begin{enumerate}
    \item $\Phi_{01}^{(1)}(x,t,k)=\begin{pmatrix}
        1\\0
    \end{pmatrix}+O(k)$ as $k\to 0$ in $\mathbb{C}^-$.

    \item $\Phi_{0 1}^{(2)}(x,t,k)=\begin{pmatrix}
        0\\1
    \end{pmatrix}+O(k)$ as $k\to 0$ in $\mathbb{C}^+$.

     \item $\Phi_{0 2}^{(1)}(x,t,k)=\begin{pmatrix}
        1\\0
    \end{pmatrix}+O(k)$ as $k\to0$ in $\mathbb{C}^+$.

    \item $\Phi_{0 2}^{(2)}(x,t,k)=\begin{pmatrix}
        0\\1
    \end{pmatrix}+O(k)$ as  $k\to0$ in $\mathbb{C}^-$.

      \item $\Phi_{0 3}^{(1)}(x,t,k)$=$\begin{pmatrix}
        1\\0
    \end{pmatrix}+O(k)$ as $k\to 0$ in $\mathbb{C}^+$.

    \item $\Phi_{03}^{(2)}(x,t,k)=\begin{pmatrix}
        0\\1
    \end{pmatrix}+O(k)$ as $k\to0$ in $\mathbb{C}^-$.

     \item $\Phi_{04}^{(1)}(x,t,k)=\begin{pmatrix}
        1\\0
    \end{pmatrix}+O(k)$ as $k\to 0$ in $\mathbb{C}^-$.

    \item $\Phi_{0 4}^{(2)}(x,t,k)=\begin{pmatrix}
        0\\1
    \end{pmatrix}+O(k)$ as $k\to 0$ in $\mathbb{C}^+$.
\end{enumerate}

Since the coefficient matrices in the Lax pair \eqref{Lax_0} are traceless, we note that
\[
\det \Phi_{0 j} \equiv 1, \qquad j=1,2,3,4.
\]
Furthermore, since $U_0(x,t,k)$, $ V_0(x,t,k)$ satisfy the symmetries \eqref{sym-UV} and $p_0(x,t,k)$ satisfy the symmetries \eqref{sym-p}, it follows that $\Phi_{0 j}(x,t,k)$ satisfy the symmetries \eqref{sym-Phi}. 

The functions $\tilde\Phi_{0 j}$,
being the solutions of both equations of \eqref{Lax_0}, are related (where they are defined) by matrices independent of $x$ and $t$, and thus

\begin{subequations}\label{rel_0}
    \begin{align}\label{rel_0_1}
        &\Phi_{0 3}(x,t,k)=\Phi_{0 2}(x,t,k)\eul^{-ikp_0(x,t,k)\sigma_3}\tilde s(k)\eul^{ikp_0(x,t,k)\sigma_3},\\\label{rel_0_2}
        &\Phi_{0 1}(x,t,k)=\Phi_{0 2}(x,t,k)\eul^{-ikp_0(x,t,k)\sigma_3}\tilde S(k)\eul^{ikp_0(x,t,k)\sigma_3},\\\label{rel_0_3}
        &\Phi_{0 4}(x,t,k)=\Phi_{0 3}(x,t,k)\eul^{-ikp_0(x,t,k)\sigma_3}\tilde S_1(k)\eul^{ikp_0(x,t,k)\sigma_3}.
    \end{align}
\end{subequations}
Due to the symmetries of $U_0$ and $V_0$,  $\tilde s(k)$ and $\tilde S(k)$
can be written as  
\begin{subequations}
    \begin{align}
        &\tilde s(k)=\Phi_{0 3}(0,0,k)=\begin{pmatrix}
        \overline{\tilde a(\bar k)} & \tilde b(k)\\
        -\overline{\tilde b(\bar k)} & \tilde a(k)
    \end{pmatrix},\\
    &\tilde S(k)=\Phi_{0 1}(0,0,k)=\begin{pmatrix}
        \overline{\tilde A(\bar k)} & \tilde B(k)\\
        -\overline{\tilde B(\bar k)} & \tilde A(k)
    \end{pmatrix},\\
    & \tilde S_1(k)=\eul^{ikL\sigma_3}\Phi_{0 4}(L,0,k)\eul^{-ikL\sigma_3}=\begin{pmatrix}
        \overline{\tilde A_1(\bar k)} & \eul^{2ikL}\tilde B_1(k)\\
        -\eul^{-2ikL}\overline{\tilde B_1(\bar k)} & \tilde A_1(k)
    \end{pmatrix}
    \end{align}
\end{subequations}

Considering  equation \eqref{inteq_03} for $t=0$ allows to derive the following properties of $\tilde a(k)$ and $\tilde b(k)$:
\begin{enumerate}

   \item $\tilde a(k)$ and $\tilde b(k)$ are analytic in $\mathbb{C}\setminus \{0\}$. Moreover, as $k\to0$
\begin{subequations}\label{s_i2}
\begin{align}
 \tilde a(k)&=1+O(k),\\
 \tilde a^*(k)&=1+O(k),\\
 \tilde b(k)&=O(k^2),\\
 \tilde b^*(k)&=O(k^2).
\end{align}
\end{subequations}     

\item The symmetries:
\begin{equation}\label{symm_atil}
    \tilde a(-k)=\overline{\tilde a(\bar k)},~ \tilde b(-k)=\overline{\tilde b(\bar k)}
\end{equation}

\item Determinant relation:

\begin{equation}\label{det_reltil}
    \tilde a(k)\overline{\tilde a(\bar k)}+ \tilde b(k)\overline{\tilde b(\bar k)}=1
\end{equation}

\end{enumerate}

Considering  equation \eqref{inteq_01} for $x=0$, we obtain the following properties of $\tilde A(k)$ and $\tilde B(k)$:

\begin{enumerate}

\item $\tilde A(k)$ and $\tilde B(k)$ are analytic in $\mathbb{C}\setminus\{0\}$. Moreover, $\tilde A(k)=1+O(k)+O\left(k\eul^{-\frac{\ii T}{2k}}
\right)$ and $\tilde B(k)=O(k)+O\left(k\eul^{-\frac{\ii T}{2k}}\right)$ as $k\to 0$.

    \item The symmetries:
    \begin{equation}\label{symm_tilA}
    \tilde A(-k)=\overline{\tilde A(\bar k)},~ \tilde B(-k)=\overline{\tilde B(\bar k)}
\end{equation}

\item Determinant relation:
\begin{equation}\label{det_rel_tilA}
    \tilde A(k)\overline{\tilde A(\bar k)}+ \tilde B(k)\overline{\tilde B(\bar k)}=1
\end{equation}

\end{enumerate}

Considering  equation \eqref{inteq_04} for $x=L$, we obtain the following properties of $\tilde A_1(k)$ and $\tilde B_1(k)$:

\begin{enumerate}
   \item $\tilde A_1(k)$ and $\tilde B_1(k)$ are analytic in $\mathbb{C}\setminus\{0\}$. Moreover, $\tilde A_1(k)=1+O(k)+O\left(k\eul^{-\frac{\ii T}{2k}}
\right)$ and $\tilde B_1(k)=O(k)+O\left(k\eul^{-\frac{\ii T}{2k}}\right)$ as $k\to 0$.

\item The symmetries:

\begin{equation}\label{symm_tilA1}
    \tilde A_1(-k)=\overline{\tilde A_1(\bar k)},~ \tilde B_1(-k)=\overline{\tilde B_1(\bar k)}
\end{equation}

\item Determinant relation:

\begin{equation}\label{det_rel_tilA1}
   \tilde  A_1(k)\overline{\tilde A_1(\bar k)}+ \tilde B_1(k)\overline{\tilde B_1(\bar k)}=1
\end{equation}

\end{enumerate}

\subsection{Relations between eigenfunctions of different Lax pairs}

Since $\Phi_0$ and $\Phi_\infty$ 
come from the 
 same system of ODEs \eqref{Lax}, they are related as follows:
\begin{subequations}\label{Phi_0_inf}
\begin{equation}\label{Phi_0_inf_}
 \Phi_{\infty j}(x,t,k)=P(x,t)\Phi_{0 j}(x,t,k)\eul^{-ikp_0(x,t,k)\sigma_3}C_j(k)\eul^{ikp(x,t,k)\sigma_3}   
\end{equation}
with $C_j(k)$  independent of $x$ and $t$ given by (since $u_x(0,0)=u_x(L,0)=0$ implies $P(0,0)=P(L,0)=I$): 
\begin{align}\label{Phi_0_inf_coeff}
    C_1&=\eul^{-\frac{\ii }{4k}T\sigma_3}P^{-1}(0,T)\eul^{-\ii k\left(\int_0^T \frac{u^2q}{2}(0,\zeta)\dd\zeta-\frac{T}{4k^2} \right)\sigma_3},\\
    C_2&=I,\\\label{Phi_0_inf_coeff3}
    C_3&=\eul^{\ii k\left(L-\theta\right)\sigma_3},\\
    C_4&=\eul^{\ii k\left(L-\frac{T }{4k^2}\right)\sigma_3}P^{-1}(L,T)\eul^{-\ii k\left(\theta-\frac{T}{4k^2}+\int_0^T \frac{u^2q}{2}(0,\zeta)\dd\zeta \right)\sigma_3}.
\end{align}
\end{subequations}

\begin{remark}\label{rem5.7}
    Substituting $(x,t)=(0,0)$ in \eqref{Phi_0_inf_} yields 
    \begin{equation}\label{s_til_s}
        s(k)=\tilde s(k)\eul^{\ii k\left(L-\theta\right)\sigma_3}.
    \end{equation}
\end{remark}

\section{Global relation }\label{sec:3}

Since the initial and boundary values $u(x,0)$, $u(0,t)$,
$u_x(0,t)$, $u(L,t)$, and
$u_x(L,t)$ cannot be prescribed as independent 
data for a well-posed problem for the SP equation on the interval, it follows 
that if these values correspond to some solution $u(x,t)$ of the SP equation in 
$x\in [0,L]$, $t\in [0,T]$, then the spectral functions
$a(k)$, $b(k)$, $A(k)$, $B(k)$, $A_1(k)$, $B_1(k)$ associated with 
them must be interrelated.

\subsection{Global relation related to $\infty$}

Expressing $\Phi_{\infty 4}(0,T,k)$ in terms of  $\Phi_{\infty 1}(0,T,k)$ using \eqref{s}--\eqref{S_1} and taking into account that
\begin{equation}\label{rel_}
\eul^{ikp(0,T,k)\hat\sigma_3}\Phi_{\infty 4}(0,T,k)=S^{-1}(k)s(k)S_1(k).
\end{equation}

Now, taking into account that $\Phi_{\infty 4}(0,T,k)$ can be 
estimated from 
    \begin{equation*}
\Phi_{\infty 4}(0,T,k)=I-
\int_{0}^{L}
	\eul^{ik\int_0^yq(\xi,T)\dd \xi\hat\sigma_3}( U_\infty\Phi_{\infty4})(y,T,k) \dd y,
\end{equation*}
we deduce, in particular, that 
$\Phi_{\infty 4}^{(12)} = O(1/k)$ as $k\to\infty$, $\Im k\ge 0$,
where $\Phi_{\infty 4}^{(12)}$ is the 
 $(12)$ entry of $\Phi_{\infty 4}(0,T,k)$.
Then, considering the $(12)$ entry of \eqref{rel_}, it follows  that 
\begin{equation}\label{gl_rel}
    B_1(k)(A(k)\overline{a(\bar k)}+B(k)\overline{b(\bar k)})\eul^{2ik\theta}+A_1(k)(A(k)b(k)-B(k)a(k))=\eul^{2ikp(0,T,k)}\ord\left(\frac{1}{k}\right), \quad\  k\to\infty,\quad \Im k> 0.
\end{equation}

\subsection{Global relation related to $0$}

Similarly, from the relations amongst $\Phi_{0j}$ we deduce that 
\begin{equation}\label{rel_i2}
    \eul^{ikp_0(0,T,k)\sigma_3}\Phi_{0 4}(0,T,k)\eul^{-ikp_0(0,T,k)\sigma_3}=\tilde S^{-1}(k)\tilde s(k)\tilde S_1(k).
\end{equation}
Here, $\Phi_{ 04}(0,T,k)$ can be 
estimated from

\begin{equation*}
\Phi_{0 4}(0,T,k)=I-
\int_{0}^{L}
	\eul^{iky\hat\sigma_3}( U_0\Phi_{04})(y,T,k) \dd y,
\end{equation*}
which gives that $\Phi_{ 04}{(12)}(0,T,k) = O(k)$ as $k\to0$
and thus

\begin{equation}\label{gl_rel_til}
    \tilde B_1(k)(\tilde A(k)\overline{\tilde a(\bar k)}+\tilde B(k)\overline{\tilde b(\bar k)})\eul^{2ikL}+\tilde A_1(k)(\tilde A(k)\tilde b(k)-\tilde B(k)\tilde a(k))=\eul^{-\frac{\ii T}{2k}}\ord\left(k\right), \quad k\to0, \quad\Im k>0.
\end{equation}

\section{Pre-Riemann-Hilbert problem  for the problem on an interval}\label{sec:4}

It is possible to combine columns of the dedicated solutions of  the Lax pair 
equations in matrices, different in different domains of the 
complex plane of the spectral parameter $k$, such that (i) 
in each domain, corresponding matrix-valued function is bounded and continuous 
up to the boundary of the domain (including unbounded domains), (ii)
their values at the boundaries of the domains are related to each other in a multiplicative manner, 
by matrices that can be determined by the initial and boundary values of the 
solution of the SP equation,  and (iii) the determinant of all matrices equals $1$
for all $x$, $t$, and $k$.

\begin{figure}[ht]
    \centering
\begin{tikzpicture}[scale=2]
\draw[red] (0,0) circle (1);
\draw[thick] (-2,0) -- (2,0);
\draw[thick, postaction={decorate},
      decoration={markings, mark=at position 0.5 with {\arrow[scale=2]{>}}}]
      (-2,0) -- (-1,0);
\draw[thick, postaction={decorate},
      decoration={markings, mark=at position 0.7 with {\arrow[scale=2]{<}}}]
      (-1,0) -- (0,0);  
\draw[thick, postaction={decorate},
      decoration={markings, mark=at position 0.5 with {\arrow[scale=2]{>}}}]
      (1,0) -- (2,0);  
\draw[thick,  postaction={decorate},
      decoration={markings, mark=at position 0.25 with {\arrow[scale=2]{<}}}] 
      (0,0) circle (1);
\draw[thick, postaction={decorate},
      decoration={markings, mark=at position 0.75 with {\arrow[scale=2]{>}}}] 
      (0,0) circle (1);
\node at (0,0.1) {\textcolor{red}{$0$}};
\fill[red] (0,0) circle(0.02);
\node at (1.1,0.1) {\textcolor{blue}{$\epsilon$}};
\fill[blue] (1,0) circle(0.02);
\end{tikzpicture}
    \caption{The oriented contour for the pre-Riemann--Hilbert problems }
    \label{fig:contour}
\end{figure}
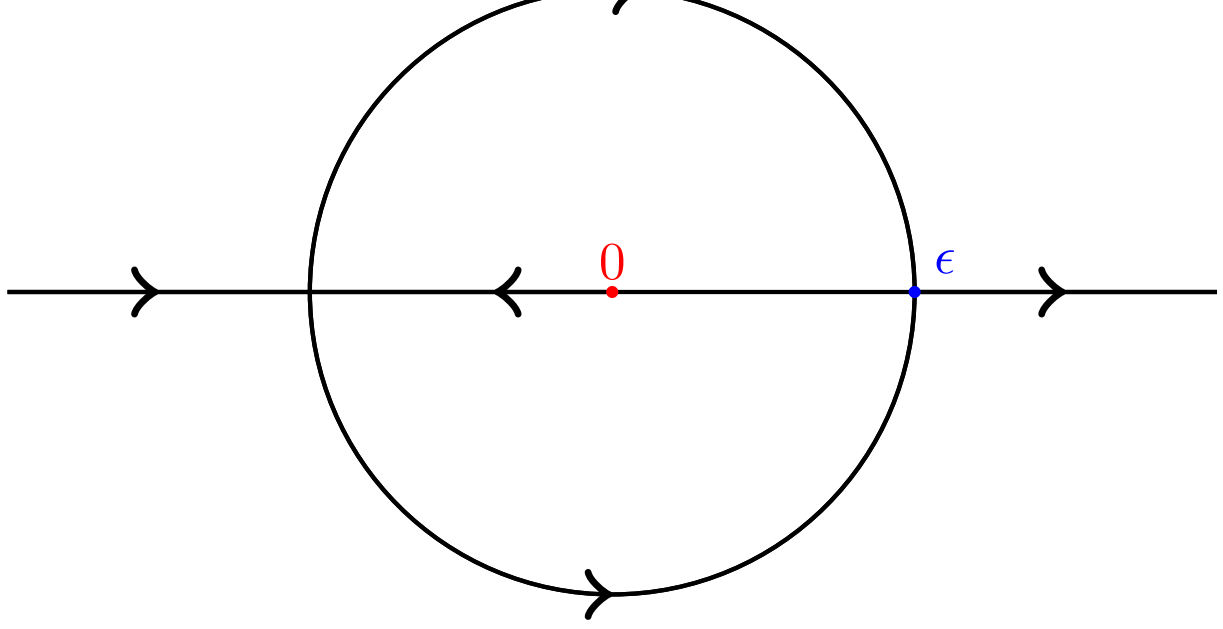

In order to achieve this introduce 
\[ 
d_1(k)=a(k)A_1(k)-b^*(k)B_1(k)\eul^{2\ii k \theta},
 \quad 
\tilde d_1(k):=\tilde a(k)\tilde A_1(k)-\tilde b^*(k)
\tilde B_1(k)\eul^{2\ii kL},
\]
and consider the following matrix-valued function defined in the domains 
in Figure \ref{fig:contour}:

\begin{equation}\label{N}
N(x,t,k)=\begin{cases}
\left( \Phi_{\infty 2}^{(1)}(k),\frac{\Phi_{\infty 4}^{(2)}(k)}{d_1(k)}\right),\quad k\in\{ \Im k > 0, |k|>\epsilon\},\\

P(x,t)\left( \Phi_{0 2}^{(1)}(k),\frac{\Phi_{0 4}^{(2)}(k)}{\tilde d_1(k)}\right)\eul^{-\ii k(p_0(x,t,k)-p(x,t,k))\sigma_3},\quad k\in\{ \Im k > 0, |k|<\epsilon\},\\

P(x,t)\left( \frac{\Phi_{0 4}^{(1)}(k)}{\tilde d^*_1(k)},\Phi_{0 2}^{(2)}(k)\right)\eul^{-\ii k(p_0(x,t,k)-p(x,t,k))\sigma_3},\quad k\in\{ \Im k < 0, |k|<\epsilon\},\\

\left(\frac{\Phi_{ \infty 4}^{(1)}(k)}{d^*_1( k)} ,\Phi_{\infty 2}^{(2)}(k)\right),\quad k\in\{ \Im k < 0, |k|>\epsilon\}.
\end{cases}
\end{equation}

\subsection{Jump conditions} By direct calculations, the boundary values of $N$ are related on $\mathbb{R}\cup\{|k|=\epsilon\}$
by 
\begin{equation}\label{jump-xt}
N_+(x,t,k)=N_-(x,t,k)J(x,t,k), \quad k\in \mathbb{R}\cup\{|k|=\epsilon\},
\end{equation}
where $J(x,t,k) = e^{-ikp(x,t,k)\sigma_3}J_0(k)e^{ikp(x,t,k)\sigma_3}$ with 
\begin{equation}\label{j0-xt}
    J_0(k)=
    \begin{cases}
\begin{pmatrix}
    1& \Gamma_2(k)\\
   0 &1
\end{pmatrix}, & \Im k>0, |k|=\epsilon,\\

\begin{pmatrix}
          1& 0\\
         \frac{\Gamma^*_1(k)}{a^{*2}(k)}&1  
\end{pmatrix}
\begin{pmatrix}
          1& r(k)\\
         r^*(k)&1  + r(k) r^*(k)
\end{pmatrix}
\begin{pmatrix}
          1& \frac{\Gamma_1(k)}{a^{2}(k)}\\
         0&1  
\end{pmatrix}
, &  \Im k = 0, |k|>\epsilon,\\
\begin{pmatrix}
          1& -\frac{\tilde\Gamma_1(k)}{\tilde a^{2}(k)}\\
         0&1  
\end{pmatrix}
\begin{pmatrix}
          1 + \tilde r(k) \tilde r^*(k)& -\tilde r(k)\\
         -\tilde r^*(k)&1  
\end{pmatrix}
\begin{pmatrix}
          1& 0\\
         -\frac{\tilde \Gamma^*_1(k)}{\tilde a^{*2}(k)}&1   
\end{pmatrix}, &  \Im k = 0, |k|<\epsilon,\\
\begin{pmatrix}
    1& 0\\
   \Gamma_2^*(k) &1
\end{pmatrix}, & \Im k<0, |k|=\epsilon.

\end{cases}
\end{equation}
Here the following notations are used:
$
r(k)=\frac{b(k)}{a(k)}
$, $\tilde r(k)=\frac{\tilde b(k)}{\tilde a(k)}$

\begin{equation}\label{Ga1}
 \Gamma_1(k)=\frac{B_1(k)a(k)\eul^{2\ii k \theta}}{a(k)A_1(k)-b^*(k)B_1(k)\eul^{2\ii k \theta}}
 =\frac{ R_1(k) a(k)\eul^{2\ii k \theta}}{ a(k)-
 b^*(k)R_1(k)\eul^{2\ii k \theta}},       
\end{equation}

\begin{equation}\label{til-Ga1}
    \tilde \Gamma_1(k)=\frac{\tilde B_1(k)\tilde a(k)\eul^{2\ii k L}}{\tilde a(k)\tilde A_1(k)-\tilde b^*(k)\tilde B_1(k)\eul^{2\ii k L}}
=\frac{\tilde R_1(k)\tilde a(k)\eul^{2\ii k L}}{\tilde a(k)-\tilde b^*(k)\tilde R_1(k)\eul^{2\ii k L}}, 
\end{equation}
\begin{equation}\label{Ga2}
    \Gamma_2 = \frac{(\tilde A_1 B_1-\tilde B_1 A_1)
\eul^{\ii k(L+\theta)}}{d_1\tilde d_1}
=\eul^{\ii k (L+\theta)}\frac{R_1-\tilde R_1}{(a-b^*R_1\eul^{2\ii k \theta})
(\tilde a-\tilde b^*\tilde R_1\eul^{2\ii k L})},
\end{equation}
where 
\begin{equation}\label{RRR}
     R_1(k):=\frac{B_1(k)}{A_1(k)},  \quad 
\tilde R_1(k)=\frac{\tilde B_1(k)}{\tilde A_1(k)}.
\end{equation}

Notice that the spectral functions associated with the boundary values 
enter the jump matrix $J_0$ through the  ratios $R_1$ and $\tilde R_1$.

\subsection{Symmetries} The symmetries \eqref{sym-Phi}, \eqref{symm_a},  and \eqref{symm_A1} imply that   $N$ satisfies the  symmetries

\begin{equation}
    \label{sym_N}
   N(k)= \overline{N(-\bar k)}=\sigma_2 \overline{N(\bar k)}\sigma_2.
\end{equation}

\subsection{Behavior at $k=\infty$} 
As $k \to \infty$, we have
\[
N(k) = I + O\!\left(\tfrac{1}{k}\right).
\]
Indeed,
from \eqref{s_inf} and \eqref{S_1_inf} it follows that
\[
d_1(k) = 1 + O\!\left(\tfrac{1}{k}\right).
\]
Together with the asymptotic properties of $\Phi_{\infty j}$, this yields the desired expansion for $N(k)$.

\subsection{Residue  conditions} 
Let $d_1(k)$ have simple zeros at $\nu_j \in \mathbb{C}^+\setminus\{|k|<\epsilon\}$. Then, by symmetry, $d^*_1( k)$ has simple zeros at $\bar\nu_j \in \mathbb{C}^-\setminus\{|k|<\epsilon\}$.

Under these assumptions, $N$ satisfies the residue conditions

\begin{align}\label{res_N-1}
\Res_{\nu_j}N^{(2)}&=\frac{\eul^{-2\ii \nu_j p(\nu_j)}}{a^2(\nu_j)}N^{(1)}(\nu_j) \Res_{\nu_j} \Gamma_1 , ~ \Res_{\nu_j}\Gamma_1=\frac{B_1(\nu_j)  \eul^{2 \ii \nu_j \theta} a(\nu_j)}{\dot d_1(\nu_j)}, ~ \nu_j\in \mathbb{C}^+, \\
\label{res_N-3}
 \Res _{\bar\nu_j}N^{(1)}&=-\frac{\eul^{2\ii \bar\nu_j p(\bar\nu_j)}}{a^{*2}(\bar\nu_j)}N^{(2)}(\bar\nu_j)  \Res _{\bar\nu_j} \Gamma_1^* ,~ \Res _{\bar\nu_j}\Gamma_1^*=\frac{B_1^*(\bar\nu_j)  \eul^{-2 \ii \bar\nu_j \theta} a^*(\bar\nu_j)}{\dot d_1^*(\bar\nu_j)},~ \bar\nu_j\in \mathbb{C}^-. 
\end{align}
Indeed,
\[
\Phi_{\infty 4}^{(2)}(x,t,k)=\Phi_{\infty 2}^{(1)}(x,t,k)\frac{B_1(k)\eul^{2\ii k\theta}}{a(k)}\eul^{-2\ii k p(x,t,k)}+\frac{d_1(k)}{a(k)}\Phi_{\infty 3}^{(2)}(x,t,k).
\]
At $k=\nu_j$ with $d_1(\nu_j)=0$, this reduces to
\[
\Phi_{\infty 4}^{(2)}(x,t,\nu_j)=\Phi_{\infty 2}^{(1)}(x,t,\nu_j)\frac{B_1(\nu_j)\eul^{2\ii \nu_j\theta}}{a(\nu_j)}\eul^{-2\ii \nu_j p(x,t,\nu_j)}.
\]
Combining this with $\Res_{\nu_j} N^{(2)} 
   = \frac{\Phi_{\infty 4}^{(2)}(\nu_j)}{\dot d_1(\nu_j)}$, \eqref{res_N-1} follows.

Similarly, one can obtain the residue conditions in the lower half-plane.

\subsection{Behavior at $k=0$}

In this subsection, $0<t<T$ is fixed and the expansions are taken non-tangentially in the indicated half-planes. Endpoint terms that are exponentially small in these limits are understood. At $t=0$, the initial values will instead be verified by the RH transformation in Section \ref{sec:6}.

Expanding $\eul^{-\ii k(p_0(x,t,k)-p(x,t,k))}$ as $k\to 0$, we obtain
\begin{equation}\label{e_i}
   \eul^{-\ii k(p_0(x,t,k)-p(x,t,k))}= 1-\ii k\left(x-\int_0^xq(\xi,t)\dd\xi-\int_0^t \frac{u^2q}{2}(0,\tau)\dd\tau\right)+O\left(k^2\right).
\end{equation}

For $k\in\mathbb{C}^+$, the properties of $\tilde s$ and $\tilde S_1$ imply 
\begin{equation}\label{d1_i2}
     \tilde d_1=1+O(k).
\end{equation}
 Expanding $\Phi_{04}^{(2)}(x,t,k)$ at $k=0$ via the Neumann series shows that
\begin{subequations}\label{beh_phi_0_0}
    \begin{equation}
    \Phi_{04}^{(2)}(x,t,k)=\begin{pmatrix}
    0\\1
\end{pmatrix}
 -\ii k \begin{pmatrix}
    u(x,t)\\0
\end{pmatrix} +O(k^2) 
    \end{equation}
Similarly, expanding $\Phi_{02}^{(1)}(x,t,k)$ at $k=0$ using the Neumann series implies
       \begin{equation}
\Phi_{02}^{(1)}(x,t,k)=\begin{pmatrix}
    1\\0
\end{pmatrix}
 -\ii k \begin{pmatrix}
    0\\ u(x,t)
\end{pmatrix} +O(k^2).
   \end{equation}
\end{subequations}

Combining \eqref{e_i}, \eqref{d1_i2} and \eqref{beh_phi_0_0}, we obtain
\begin{equation}\label{Psi_i2_1}
N(x,t,k)=P(x,t)\left(I-\ii k\begin{pmatrix}
     \tilde \gamma_1(x,t)&u(x,t)\\
     u(x,t)& \tilde \gamma_2(x,t)
\end{pmatrix}+O(k^2)\right).
\end{equation}

Similarly, for $k\in\mathbb{C}^-$, \eqref{Psi_i2_1} holds.

\section{The SP equation,  associated Lax pair, and  RH formalism in $y,t$ variables}\label{sec:5}

\subsection{The SP equation and the  associated Lax pair in $y,t$ variables}

In order to have the data for the RH problem to 
be  explicitly determined by the initial data only,  
it is convenient to introduce the new space variable $y=y(x,t)$ in such a way that  
\begin{subequations}\label{y-xt}
\begin{align}
\label{y_x} y_x & = q,\\
   \label{y_t} y_t &=\frac{u^2q}{2}
\end{align}
    \end{subequations}
    which is compatible due to \eqref{SP}. These relations determine $y$ up to an additive constant, which is fixed differently depending on the spatial domain.
    Particularly, for the problem on the interval $x\in [0,L]$ we set
    \begin{equation}\label{y_}
y(x,t)=\int_0^xq (\xi,t)\dd\xi+\int_0^t\frac{u^2q}{2}(0,\tau)\dd \tau,
\end{equation}
which is fixed by the condition $y(0,0)=0$.
 
Then, by differentiating the identity $x(y(x,t),t)=x$, we derive the equations that describe the 
inverse mapping $x=x(y,t)$:
\begin{subequations}\label{x-yt}
\begin{align}
\label{x_y} x_y & = \frac{1}{\hat q},\\
   \label{x_t} x_t &=-\frac{\hat u^2}{2},
\end{align}
    \end{subequations}
where $\hat u (y,t) = u(x(y,t),t)$ and  $\hat q (y,t) = q(x(y,t),t)$.

Notice that for the problem on $[0,L]$ with the periodic conditions, one has
\begin{equation}\label{y_0L}
    y(L,t)= y(0,t)+\theta\text{ for all } t,
\end{equation}
where
$\theta= \int_0^L q(\xi,0)\dd\xi $.

\begin{proposition} (SP equation in the $(y,t)$ variables) Let $u(x,t)$, $w(x,t)$  and $q(x,t)$ satisfy \eqref{SP} and let $y(x,t)$ be 
such that \eqref{y-xt} hold.
Then the SP equation \eqref{SP} in the $(y,t)$ variables reads as the following system of equations
for  $\hat q(y,t):=q(x(y,t),t)$, $\hat u(y,t):=u(x(y,t),t)$, and 
$\hat w(y,t):=w(x(y,t),t)$:

\begin{subequations}\label{SP_in_y}
\begin{align}\label{SP_in_y-1}
\hat q_t&=\hat u \hat q\hat w,\\\label{SP_in_y-2}
\hat q&=\sqrt{\hat w^2+1},\\\label{SP_in_y-3}
\hat u_y&=\frac{\hat w}{\hat q},\\
\label{SP_in_y-4}
\hat w_t&=\hat u\hat q^2.
\end{align}
\end{subequations}
    
\end{proposition}

\begin{proof}
As we discussed above, \eqref{y-xt} implies \eqref{x-yt}. 
 Substituting from \eqref{SP-2} and \eqref{x_t} into the equality
\[
\hat q_t=\left(q_xx_t+q_t\right)\mid_{x=x(y,t)}
\]
we get \eqref{SP_in_y-1}.

Now, substituting \eqref{x_y} into the equality
$\hat u_y(y,t)=\left(u_x(x,t)x_y(y,t)\right)\mid_{x=x(y,t)}$
we get \eqref{SP_in_y-3}.

Finally, $\hat q(y,t)=\sqrt{w^2(x(y,t),t)+1}=\sqrt{\hat w^2(y,t)+1}$ and thus \eqref{SP_in_y-2} follows. Moreover, \eqref{SP-w-evolution} and \eqref{x_t} imply
\[
    \hat w_t=\left(w_t-\frac{u^2}{2}w_x\right)\big|_{x=x(y,t)}
    =\hat u\hat q^2,
\]
which gives \eqref{SP_in_y-4}.
\end{proof}

\begin{remark}\label{rem:cons_law}
Notice that we can rewrite \eqref{SP_in_y-1}--\eqref{SP_in_y-3} as
\begin{equation}\label{SP_in_y_alt}
    \left(\frac{1}{\hat q}\right)_t=-\frac{1}{2}\left(\hat u^2\right)_y, \quad \hat q=\sqrt{\hat w^2+1}, \quad \hat u_y=\frac{\hat w}{\hat q},
\end{equation}
which implies that for the periodic (with period $\theta$) solution the following conservation law holds:
\begin{equation*}
    \int_0^\theta \frac{1}{\hat q}(\xi,t)\dd \xi=\int_0^\theta \frac{1}{\hat q}(\xi,0)\dd \xi.
\end{equation*}
\end{remark}

Similarly, the reverse  change of variable $(y,t)\mapsto(x,t)$ 
reduces \eqref{SP_in_y} to \eqref{SP} with 
$u(x,t):=\hat u(y(x,t),t)$, $w(x,t):=\hat w(y(x,t),t)$ and $q(x,t):=\hat q(y(x,t),t)$. 
\begin{proposition}\label{prop:hatu_u}
    Let $\hat u(y,t)$, $\hat w(y,t)$,  and $\hat q(y,t)$ satisfy 
    \eqref{SP_in_y} and let $x(y,t)$ be 
such that \eqref{x-yt} holds. Define $u(x,t):=\hat u(y(x,t),t)$, $w(x,t):=\hat w(y(x,t),t)$ and $q(x,t):=\hat q(y(x,t),t)$.
Then \eqref{SP} holds for $u(x,t)$, $w(x,t)$ and $q(x,t)$.
\end{proposition}

Now let us show that a solution of the periodic problem for the SP equation 
in the $x,t$ variables \eqref{SP} gives rise to a periodic solution 
of the SP equation 
in the $y,t$ variables \eqref{SP_in_y} and vice versa.

\begin{proposition}\label{prop:periodic1}
    Let  $u(x,t)$, $w(x,t)$ and $q(x,t)$  be the components of a periodic solution of  \eqref{SP} with the period $L$, 
        and let $y(x,t)$ be defined by \eqref{y_}. Then 
    $\hat u(y,t):=u(x(y,t),t)$, $\hat w(y,t):=w(x(y,t),t)$, and 
$\hat q(y,t):=q(x(y,t),t)$ are the components of a periodic solution of 
\eqref{SP_in_y}, with the period
\begin{equation}\label{theta}
    \theta:= \int_0^L q(\xi,0) d\xi.
\end{equation}
\end{proposition}

\begin{proof}
    Introducing $x(y,t)$ as the inverse to $y(x,t)$ and setting $y=0$ and $y=\theta$
    in the l.h.s. of \eqref{y_} we have
    \begin{align*}
        0 & = \int_0^{x(0,t)} q(\xi,t) d\xi 
            +\int_0^t\frac{u^2q}{2}(0,\tau)\dd \tau, \\
            \theta & = \int_0^{x(\theta,t)}  q(\xi,t) d\xi 
            +\int_0^t\frac{u^2q}{2}(0,\tau)\dd \tau,
    \end{align*}
    from which we deduce that 
    \[
    \int_{x(0,t)}^{x(\theta,t)} \sqrt{1+w^2(\xi,t)} d\xi =\theta = \int_0^L \sqrt{1+w^2(\xi,0)} d\xi
    \]
    for all $t$. Since $q(x,t)$ is periodic and $q(x,t)>0$, it follows that 
    \[
    x(\theta,t) - x(0,t) = L \ \ \text{for all} \ \ t.
    \]
    Consequently,
    \[
    \hat u(\theta,t) =  u(x(\theta,t),t) = u(x(0,t),t) = \hat u(0,t)
    \]
    and thus $\hat u(y,t)$ is periodic as function of $y$, for all fixed $t$, with period $\theta$. Similarly for $\hat w$ and $\hat q$.
\end{proof}

\begin{proposition}\label{prop:periodic2}
    Let  $\hat u(y,t)$, $\hat w(y,t)$,  and $\hat q(y,t)$ be the components of a periodic solution of  \eqref{SP_in_y} with the period $\theta$, 
        and let $x(y,t)$ be defined by 
\begin{equation}
    \label{x}
    x(y,t)=\int_0^y \frac{1}{\hat q(\xi,t)}\dd\xi-
    \int_0^t \frac{\hat u^2(0,\zeta)}{2}\dd\zeta.
\end{equation}
Then 
    $u(x,t):=\hat u(y(x,t),t)$, $w(x,t):=\hat w(y(x,t),t)$ and $q(x,t):=\hat q(y(x,t),t)$
    are the components of a periodic solution of 
\eqref{SP}, with the period
\begin{equation}\label{L}
    L= \int_0^\theta \frac{1}{\hat q(\xi,0)}\dd\xi.
\end{equation}
    Moreover, \eqref{theta} holds true.
\end{proposition}

\begin{proof}
    Introducing $y(x,t)$ as the inverse to $x(y,t)$ and arguing as in the proof of Proposition
    \ref{prop:periodic1} we arrive at
    \[
    \int_{y(0,t)}^{y(L,t)} \frac{1}{\hat q(\xi,t)} d\xi =L
    = \int_0^\theta \frac{1}{\hat q(\xi,0)}  d\xi
    \]
    for all $t$, which implies  
    \[
    y(L,t) - y(0,t) = \theta \ \ \text{for all} \ \ t.
    \]
     Consequently,
    \[
    u(L,t) =  \hat u(y(L,t),t) = \hat u(y(0,t),t) = u(0,t)
    \]
    and thus $u(x,t)$ is periodic as function of $x$, for all fixed $t$, with period $L$
    determined by \eqref{L}.  Similarly for $w$ and $q$.

    Finally,
    \[
    \int_0^L q(\xi,0) d\xi = \int_0^L \hat q(y(\xi,0),0) d\xi=
    \int_0^\theta \hat q(y,0) x_y(y,0)\dd y = \int_0^\theta dy =\theta.
    \]
\end{proof}

Now, let us reformulate the original Lax pair equations in the $(y,t)$ variables.

\begin{proposition}
The Lax pair \eqref{Lax-Q-form} in the variables $(y,t)$ takes the form 
\begin{subequations}\label{Lax_y}
\begin{align}\label{Lax_y_y}
    &\hat \Psi_{ y}=\left(-\ii k \sigma_3+\frac{\hat w_y}{2\hat q^2}\begin{pmatrix}
        0&1\\-1&0
    \end{pmatrix}\right)\hat \Psi,\\\label{Lax_y_t}
   &\hat\Psi_{ t}=-\frac{1}{4\ii k}\begin{pmatrix}
        \frac{1}{\hat q}&-\frac{\hat w}{\hat q}\\-\frac{\hat w}{\hat q}&-\frac{1}{\hat q}
    \end{pmatrix} \hat\Psi.
\end{align}
\end{subequations}
\end{proposition}

\begin{proof}
  Introducing $\hat\Psi(y,t) = \tilde\Phi_\infty (x(y,t),t)$ and taking into account \eqref{x_y} and \eqref{x_t}, the Lax pair \eqref{Lax-Q-form} in the variables $(y,t)$ takes the form \eqref{Lax_y}.
\end{proof}

\subsection{The  RH formalism in $y,t$ variables}

In this subsection, we describe the derivation of a (local) solution to the SP equation in the $(y,t)$ variables, starting from an RH problem parametrized by $y$ and $t$, motivated by the considerations in Section \ref{sec:4}.

Let $\Gamma$ denote an oriented contour in the complex plane, invariant under the symmetries $k\mapsto-k$ and $k\mapsto\bar k$. The plus and minus sides of each arc are its left and right sides, respectively. Orient the symmetric arcs so that $k\mapsto-\bar k$ preserves the plus and minus sides, while $k\mapsto\bar k$ interchanges them. This local convention also applies to finite branch cuts.

Consider the following \textbf{RH problem parametrized by $y$ and $t$}: Given the $2\times 2$-matrix valued function $J_0(k)$ for $k\in\Gamma$ such that 
\begin{subequations}\label{jump_cond}
    \begin{align}\label{det_jump}
    &\det J_0(k)\equiv 1,\\
    \label{jump_at_inf}
    &J_0(k)=I+O(\frac{1}{k}),\quad k\to\infty,\\
    \label{jump_at_0}
    &J_0(k)=I+O(k),\quad k\to0,
    \end{align}
\end{subequations}
and the sets $\{k_j,c_j\}_1^N$ with $k_j\in\mathbb C\setminus\Gamma$, and $\{\eta_j,d_j\}_1^M$ with $\eta_j\in\mathbb C\setminus\Gamma$,
find a piece-wise meromorphic (w.r.t. $\Gamma$) $2\times 2$-matrix valued function $\hat M(y,t,k)$ satisfying:
\begin{enumerate}[\textbullet]
\item
\emph{Jump} condition
\begin{equation}\label{jump-y_loc}
\hat M_+(y,t,k)=\hat M_-(y,t,k) \hat J(y,t,k),\qquad k\in\Gamma,
\end{equation}
where  
$\hat J(y,t,k)=\eul^{-ik\hat p(y,t,k)\sigma_3}J_0(k)\eul^{ik\hat p(y,t,k)\sigma_3}$ with  
$\hat p(y,t,k):=y-\frac{t}{4k^2}$.

\item  \emph{Normalization} condition:
\begin{equation}\label{norm-m-hat_loc}
\hat M(y,t,k)=I+\ord(\frac{1}{k}), \quad k\to\infty.
\end{equation}

\item \emph{Residue} conditions:
\begin{align}\label{res_hatM_loc}
\Res_{ k_j}\hat M^{(1)}(y,t,k)&= c_j \hat M^{(2)}(y,t,k_j)\eul^{2\ii k_j\hat p(y,t,k_j)}, ~  k_j\in\mathbb C\setminus\Gamma, \\\label{res_hatM__loc}
\Res_{\eta_j}\hat M^{(2)}(y,t,k)&=d_j \hat M^{(1)}(y,t,\eta_j)\eul^{-2\ii\eta_j\hat p(y,t,\eta_j)} , ~ \eta_j\in\mathbb C\setminus\Gamma.
\end{align}
\end{enumerate}

The poles are assumed simple, finite in number, nonzero, and disjoint from the contour. We require boundedness at $k=0$ and at any other finite contour vertices, and the usual continuous boundary values away from the vertices. In the reconstruction arguments below, we additionally assume sufficient differentiability in $y,t$, boundedness of $\hat M_y$ and $\hat M_t$ near $k=0$, and the non-tangential asymptotic expansions used below, differentiable term by term. These regularity conditions are part of the reconstruction hypotheses; they do not follow from \eqref{jump_cond} alone.

\begin{proposition}
 Assume that $\hat M $ is a solution of RH problem \eqref{jump-y_loc}--\eqref{res_hatM__loc}. Then $\det \hat M=1$.
\end{proposition}
\begin{proof}
The fact that $\det J_0(k)\equiv 1$ implies that $\det \hat M$ has neither a jump across $\Gamma$, nor singularities at the points $k_j$ and $\eta_j$, and thus is an entire function. Therefore, \eqref{jump_at_inf}
together with the Liouville's theorem implies $\det \hat M=1$.
\end{proof}

\begin{proposition}
     If a solution of the RH problem \eqref{jump-y_loc}--\eqref{res_hatM__loc} exists, it is unique.
\end{proposition}

\begin{proof}
    The statement follows by considering two solutions $\hat M_1$ and $\hat M_2$ and applying Liouville's theorem to their ratio $\hat M_1 (\hat M_2)^{-1}$.
\end{proof}

Now, assume that the RH problem \eqref{jump-y_loc}--\eqref{res_hatM__loc} 
has a solution $\hat M(y,t,k)$ that satisfies the \emph{symmetries}
\begin{equation}
 \label{sym_1}
 \hat M(k)=\overline{\hat M(-\bar k)}
 \end{equation}
and 
 \begin{equation}
 \label{sym_2}
\hat M(k)=\begin{pmatrix}
    0&-\ii\\\ii&0
\end{pmatrix}\overline{\hat M(\bar k)}\begin{pmatrix}
    0&-\ii\\\ii&0
\end{pmatrix}
 \end{equation} 
and  is differentiable w.r.t. $y$ and $t$.

\begin{remark}
The symmetries and $\det\hat M=1$ do not alone imply solvability. To invoke Zhou's vanishing lemma \cite{Z89}, one must also verify the appropriate Schwarz reflection condition off the real axis, positivity of the Hermitian part of the jump on the real axis, and the contour and function-space hypotheses. For a meromorphic problem, the poles must first be replaced by jumps on small paired circles, or treated by an equivalent meromorphic version of the lemma. These conditions are verified for the master jump below under the stated nondegeneracy assumptions.
\end{remark}

Now, our goal is to show that, by evaluating the solution $\hat M(y,t,k)$ of the RH problem \eqref{jump-y_loc}--\eqref{res_hatM__loc} at appropriate points of the extended complex plane $\overline{\mathbb C}$, one can 
    obtain a solution of the SP equation in the $y,t$ variables \eqref{SP_in_y}.
For this purpose,
starting from $\hat M(y,t,k)$, we define $2\times 2$-matrix valued functions \[ \hat \Psi(y,t,k):=\hat M(y,t,k)\eul^{-\ii k\hat p(y,t,k)\sigma_3}\] and show that $\hat\Psi(y,t,k)$ satisfies the system of differential equations:
\begin{equation}\label{Lax-hat-hat}
\begin{split}
\hat\Psi_y&=\doublehat{U}\hat\Psi, \\
\hat \Psi_t&=\doublehat{V}\hat\Psi,
\end{split}
\end{equation}
where $\doublehat{U}$ and $\doublehat{V}$ have the same (rational) dependence on $k$ as in \eqref{Lax_y_y} and \eqref{Lax_y_t}, with coefficients given in terms of $\hat M(y,t,k)$ evaluated at appropriate values of $k$.

\begin{proposition}\label{Prop_Lax_y}
    Let $\hat M(y,t,k)$ be the solution of the RH problem \eqref{jump-y_loc}--\eqref{res_hatM__loc} that satisfies symmetries \eqref{sym_1} and \eqref{sym_2}. Define $\hat \Psi (y,t,k):= \hat M(y,t,k)\eul^{-\ii k \hat p(y,t,k)\sigma_3}$. 

    Then $\hat\Psi(y,t,k)$ satisfies the differential equation
\begin{equation}\label{Lax_y_Psi}
    \hat\Psi_y=\doublehat{U}\hat\Psi
\end{equation}
with 
\begin{equation}\label{hat_hat_U}
    \doublehat{U}(y,t,k)=-\ii k \sigma_3 + 2\hat\xi(y,t)\begin{pmatrix}
    0&-1\\1&0
\end{pmatrix},
\end{equation}
where 
$\hat \xi(y,t)\in\mathbb{R}$ can be obtained from the large $k$ expansion of $\hat M(y,t,k)$: 
\begin{equation}\label{xi}
\hat \xi(y,t)=\hat m_2^\infty(y,t),
\end{equation}
where
\begin{equation}\label{xi-M}
\hat M(y,t,k)=I+\frac{\ii}{k}\begin{pmatrix}
    \hat m_1^\infty & \hat m_2^\infty\\
    \hat m_2^\infty& -\hat m_1^\infty
\end{pmatrix}(y,t)+O(\frac{1}{k^2}),\qquad k\to\infty.
\end{equation}

\end{proposition}

\begin{proof}
     First, notice that $\hat\Psi(y,t,k)$ satisfies the jump condition
\[
\hat\Psi^+(y,t,k)=\hat\Psi^-(y,t,k)J_0(k)
\]
with the jump matrix $J_0(k)$ independent of $y$. Hence, $\hat\Psi_y(y,t,k)$ satisfies the same jump condition.

As for the residue conditions, we notice that the symmetry assumptions  
\eqref{sym_1} and \eqref{sym_2} imply that if $ k_j$ is pole of $\hat M^{(1)}$ then $-\bar k_j$ is also pole of $\hat M^{(1)}$ with $c_{ k_j}=-\overline{c_{-\bar k_j}}$. Similarly, if $\eta_j$ is pole of $\hat M^{(2)}$ then $-\bar\eta_j$ is also pole of $\hat M^{(2)}$ with $d_{\eta_j}=-\overline{d_{-\bar\eta_j}}$. Moreover,
\begin{equation}\label{res-sym}
    \eta_j=- k_j, \quad d_j=c_j.
\end{equation}
This implies, $d_{\bar  k_j}=-\overline{c_{ k_j}}$.

Since $\{c_j\}$ are independent of $y$, 
$\hat\Psi_y(y,t,k)$ satisfies the same residue conditions as $\hat\Psi(y,t,k)$ does:

\begin{align*}
\Res_{ k_j}\hat \Psi^{(1)}&= c_j \hat \Psi^{(2)}( k_j) , ~  k_j\in\mathbb C\setminus\Gamma, \\
\Res_{\bar k_j}\hat \Psi^{(2)}&=-\bar{ c}_j \hat \Psi^{(1)}(\bar k_j) , ~ \bar k_j\in\mathbb C\setminus\Gamma.
\end{align*}
with constants $c_j$ independent of $y$.
Consequently,  $\hat \Psi_y \hat \Psi^{-1}=\hat M_y \hat M^{-1}-\ii k \hat p_y\hat M \sigma_3 \hat M ^{-1}$ (with $\hat p_y=1$)  has neither jumps nor singularities at $ k_j$ and thus it is a meromorphic function, with possible singularity at $k=\infty$.

Let us analyze the behavior of $\hat \Psi_y \hat \Psi^{-1}$ as $k\to\infty$.

\begin{enumerate}[(i)]
    \item As $k\to\infty$, we have due to the symmetry \eqref{sym_1}

    \begin{equation}
   \hat M = I+\frac{\ii}{k}\begin{pmatrix}
    \hat m_{11}^\infty & \hat m_{12}^\infty\\
    \hat m_{21}^\infty& \hat m_{22}^\infty
\end{pmatrix}+O(\frac{1}{k^2}), \quad k\in\mathbb{C}^+
   \end{equation}
with $\hat m_{ij}^\infty \in\mathbb{R}$. Moreover, the symmetry \eqref{sym_2} implies
    \begin{equation}
   \hat M = I-\frac{\ii}{k}\begin{pmatrix}
    \hat m_{22}^\infty & -\hat m_{21}^\infty\\
    -\hat m_{12}^\infty& \hat m_{11}^\infty
\end{pmatrix}+O(\frac{1}{k^2}), \quad k\in\mathbb{C}^-.
   \end{equation}
Consequently, $
\hat M_y=O(\frac{1}{k})$  and thus 
 $\hat M_y \hat M^{-1}=O(\frac{1}{k})$, which leads to 
the following expansion for $\hat\Psi_y\hat\Psi^{-1}$:

\begin{equation}\label{Psi_y_Psi_inf_C+}
  \hat\Psi_y\hat\Psi^{-1}
  =-\ii k\sigma_3+
  \begin{pmatrix}
  0&-2\hat m_{12}^\infty\\
  2\hat m_{21}^\infty&0
  \end{pmatrix}
  +O(k^{-1}),\qquad k\in\mathbb C^+.
\end{equation}
In $\mathbb C^-$, the symmetry gives
\begin{equation}\label{Psi_y_Psi_inf_C-}
  \hat\Psi_y\hat\Psi^{-1}
  =-\ii k\sigma_3+
  \begin{pmatrix}
  0&-2\hat m_{21}^\infty\\
  2\hat m_{12}^\infty&0
  \end{pmatrix}
  +O(k^{-1}).
\end{equation}
Since $\hat\Psi_y\hat\Psi^{-1}$ is meromorphic across the contour, these expansions coincide. Thus $\hat m_{12}^\infty=\hat m_{21}^\infty$. Also, $\det\hat M=1$ implies $\hat m_{11}^\infty=-\hat m_{22}^\infty$.

Denoting $\hat m_1^\infty=\hat m_{11}^\infty$ and $\hat m_2^\infty=\hat m_{12}^\infty$, we obtain
\begin{equation}\label{Psi_y_Psi_inf}
  \hat\Psi_y\hat\Psi^{-1}
  =-\ii k\sigma_3+
  2\hat m_2^\infty
  \begin{pmatrix}0&-1\\1&0\end{pmatrix}
  +O(k^{-1}).
\end{equation}
By Liouville's theorem, this implies \eqref{Lax_y_Psi}--\eqref{xi-M}.
\end{enumerate}
\end{proof}

\begin{remark}\label{rem:expM_1}
Under the local regularity assumptions stated above, the symmetries \eqref{sym_1} and \eqref{sym_2} lead to the following non-tangential expansions near $k=0$:
\begin{align}\label{hatM_at_0_pl}
    &\hat M=B\left(I+\ii k\Gamma_++O(k^2)\right),\qquad k\in\mathbb C^+,\\
\label{hatM_at_0_min}
    &\hat M=B\left(I-\ii k\Gamma_-+O(k^2)\right),\qquad k\in\mathbb C^-,
\end{align}
where
\[
 B=\begin{pmatrix}\hat\alpha&\hat\beta\\-\hat\beta&\hat\alpha\end{pmatrix},
 \quad
 \Gamma_+=\begin{pmatrix}
 \hat\gamma_{11}&\hat\gamma_{12}\\
 \hat\gamma_{21}&\hat\gamma_{22}
 \end{pmatrix},
 \quad
 \Gamma_-=\begin{pmatrix}
 \hat\gamma_{22}&-\hat\gamma_{21}\\
 -\hat\gamma_{12}&\hat\gamma_{11}
 \end{pmatrix},
\]
with real coefficients. The identity $\det\hat M=1$ gives
$\hat\alpha^2+\hat\beta^2=1$ and $\operatorname{tr}\Gamma_\pm=0$.
For a $2\times2$ matrix $A$, denote its adjugate by $A^\sharp$.
Thus $\Gamma_\pm^\sharp=-\Gamma_\pm$.

Substituting these expansions into
$\hat\Psi_y\hat\Psi^{-1}=\hat M_y\hat M^{-1}-\ii k\hat M\sigma_3\hat M^{-1}$
and comparing with \eqref{hat_hat_U}, we obtain
\begin{align}
 B_yB^{-1}&=2\hat\xi\begin{pmatrix}0&-1\\1&0\end{pmatrix},\\
\label{BGamma_rel}
 \pm B\Gamma_{\pm y}B^{-1}-B\sigma_3B^{-1}&=-\sigma_3.
\end{align}
\end{remark}

\begin{proposition}\label{Prop_Lax_t}
$\hat\Psi(y,t,k)$ defined as in Proposition \ref{Prop_Lax_y} satisfies the differential equation
\begin{equation}\label{Lax_t_Psi}
    \hat\Psi_t=\doublehat{V}\hat\Psi
\end{equation}
with 
\begin{equation}\label{hat_hat_V}
    \doublehat{V}=-\frac{1}{4\ii k} \begin{pmatrix}
        \hat \alpha^2-\hat \beta^2&-2\hat \alpha \hat \beta \\
        -2\hat \alpha \hat \beta&-\hat \alpha^2+\hat \beta^2
    \end{pmatrix},
\end{equation}
where 
\begin{enumerate}[(i)]
    \item 
$\hat \alpha(y,t)\in\mathbb{R}$ and $\hat \beta(y,t)\in\mathbb{R}$ can be obtained from the expansion of $\hat M(y,t,k)$ as $k\to 0$: 
\begin{equation}\label{alpha_beta-M}
    \hat M(y,t,k)=\begin{pmatrix}
    \hat \alpha &\hat \beta\\
   -\hat \beta&\hat \alpha
\end{pmatrix}+
O(k),\quad k\to 0.
\end{equation}
\end{enumerate}

\end{proposition}

\begin{proof}
    Similarly to Proposition \ref{Prop_Lax_y}, we notice that $\hat\Psi_t \hat\Psi^{-1}=\hat M_t\hat M^{-1}-\ii k\hat p_t\hat M\sigma_3\hat M^{-1}$ has neither jump nor poles at $ k_j$, and thus it is a meromorphic function, with possible singularities at $k=\infty$ and $k=0$, the latter being due to the singularity of $\hat p_t$ at $k=0$:

Evaluating $\hat\Psi_t\hat\Psi^{-1}$ near these points, we have the following.
\begin{enumerate}[(i)]

\item As $k\to\infty$, we have $\hat M_t\hat M^{-1}=O(\frac{1}{k})$, $\hat p_t(k)=-1/(4k^2)$ and thus
\begin{equation}\label{psi-inf-t}
\hat\Psi_t\hat\Psi^{-1}(k)=\ord(k^{-1}),\qquad k\to\infty.
\end{equation}

\item As $k\to 0$, using \eqref{hatM_at_0_pl} and \eqref{hatM_at_0_min}, we have

\begin{align}\label{psi_t_psi_0}
    \hat\Psi_t\hat\Psi^{-1}=\frac{\ii}{4 k}B\sigma_3B^{-1}+B_tB^{-1}-\frac{1}{4}B\left(\pm \sigma_3\Gamma_{\pm}^{\sharp}\pm\Gamma_{\pm}\sigma_3\right)B^{-1}+O(k),\quad k\to 0, \quad k\in\mathbb{C}^\pm
\end{align}
In particular,
\begin{equation}\label{Psi_t_Psi_0}
  \hat \Psi_t \hat \Psi^{-1}=-\frac{1}{4\ii k} \begin{pmatrix}
\hat \alpha^2- \hat \beta^2&-2 \hat \alpha \hat \beta\\-2 \hat \alpha \hat \beta&-\hat \alpha^2+\hat \beta^2
\end{pmatrix}+O(1),  \quad k\to 0.
\end{equation}

Now, \eqref{psi-inf-t} and \eqref{Psi_t_Psi_0} together with the Liouville's theorem imply \eqref{Lax_t_Psi}--\eqref{alpha_beta-M}.

\end{enumerate}

\end{proof}

\begin{remark}\label{rem:expM_2}
The constant term in \eqref{psi_t_psi_0} must vanish, since \eqref{hat_hat_V} has only a $k^{-1}$ term. Equivalently,
\begin{equation}\label{BGamma_rel_t}
 B^{-1}B_t=\frac14[\Gamma_+,\sigma_3].
\end{equation}
The left-hand side is real and skew-symmetric. Therefore
$\hat\gamma_{12}=\hat\gamma_{21}$.
Together with $\operatorname{tr}\Gamma_+=0$, this gives
\begin{equation}\label{Gamma_pm}
 \Gamma_+=\Gamma=-\Gamma_-=
 \begin{pmatrix}\hat f_1&\hat f_2\\\hat f_2&-\hat f_1\end{pmatrix},
 \qquad \Gamma^\sharp=-\Gamma.
\end{equation}
Consequently,
\begin{equation}\label{M_at_0_exp}
 \hat M(y,t,k)=
 \begin{pmatrix}\hat\alpha&\hat\beta\\-\hat\beta&\hat\alpha\end{pmatrix}
 \left(I+\ii k
 \begin{pmatrix}\hat f_1&\hat f_2\\\hat f_2&-\hat f_1\end{pmatrix}
 +O(k^2)\right),\qquad k\to0.
\end{equation}
\end{remark}

\begin{proposition}
    Let $\hat \alpha(y,t)$, $\hat \beta(y,t)$, $\hat f_1(y,t)$ and $\hat f_2(y,t)$ be the functions determined in terms of $\hat M(y,t,k)$ as in Propositions \ref{Prop_Lax_y} and \ref{Prop_Lax_t}. Then they satisfy the following equations:
    \begin{subequations}\label{exp_rel}
        \begin{align}
            \label{exp_rel_1}
&\hat f_{1y}=1-(\hat\alpha^2-\hat \beta^2),\\\label{exp_rel_2}
&\hat f_{2y}=-2\hat\alpha\hat \beta,\\\label{exp_rel_4}
&\hat \alpha_{t}=\frac{\hat f_2}{2}\hat \beta,\\\label{exp_rel_5}
&\hat \beta_{t}=-\frac{\hat f_2}{2}\hat \alpha,\\
\label{exp_rel_3}
&\hat \alpha^2+\hat \beta^2=1.
        \end{align}
    \end{subequations}
\end{proposition}

\begin{proof}
First, notice that \eqref{exp_rel_3}  follows directly from $\det \hat M=1$.
Next, substituting \eqref{Gamma_pm} into \eqref{BGamma_rel} yields \eqref{exp_rel_1} and \eqref{exp_rel_2} while substituting \eqref{Gamma_pm} into \eqref{BGamma_rel_t} yields \eqref{exp_rel_4} and \eqref{exp_rel_5}
\end{proof}

\begin{proposition}
    \label{reduce} If $\hat\alpha^2-\hat\beta^2>0$ (then \eqref{exp_rel_3} implies that $0<\hat\alpha^2-\hat\beta^2\leq1$).
Introducing  $\hat u(y,t)$, $\hat w(y,t)$, and 
$\hat q(y,t)$
in terms of $\hat \alpha$, $\hat \beta$,  $\hat f_1$, and $\hat f_2$ by
\begin{equation}\label{hmbbgg}
\hat u=-\hat f_2,\quad \hat w = \frac{2\hat\alpha\hat\beta}{\hat\alpha^2-\hat\beta^2},\quad \hat q = \frac{1}{\hat\alpha^2-\hat\beta^2},
\end{equation}
 equations  \eqref{exp_rel} -- \eqref{exp_rel_3} reduce to \eqref{SP_in_y}.
\end{proposition}
    
\begin{proof}
    First, notice that \eqref{exp_rel_4}, \eqref{exp_rel_5}, and \eqref{exp_rel_3} imply
\[
    \hat f_2=2(\hat\beta\hat\alpha_t-\hat\alpha\hat\beta_t),
\]
which does not require division by $\hat\alpha$ or $\hat\beta$.
    Now, differentiating  $\hat q$ w.r.t. $t$, we get
    \[\hat q_t=\frac{-2(\hat \alpha\hat \alpha_t-\hat \beta\hat \beta_t)}{(\hat \alpha^2-\hat \beta^2)^2}=\frac{-2\hat \alpha\hat \beta\hat f_2}{(\hat \alpha^2-\hat \beta^2)^2}=\hat u\hat w\hat q,\]
    which gives \eqref{SP_in_y-1}.

    Now, let us compute $\hat u_y$ using \eqref{exp_rel}:
    \[\hat u_y=-\hat f_{2y}=2\hat\alpha\hat\beta=\frac{\hat w}{\hat q}\]
    which gives \eqref{SP_in_y-3}.

    Moreover, writing $n=2\hat\alpha\hat\beta$ and $d=\hat\alpha^2-\hat\beta^2$, we have $n_t=\hat u d$, $d_t=-\hat u n$, and $n^2+d^2=1$. Hence
    \[
       \hat w_t=\left(\frac nd\right)_t
       =\frac{\hat u(d^2+n^2)}{d^2}
       =\hat u\hat q^2,
    \]
    which gives \eqref{SP_in_y-4}.

    In order to get \eqref{SP_in_y-2} notice that due to \eqref{exp_rel_3} we have $1+\hat w^2=1+\frac{4\hat\alpha^2\hat\beta^2}{(\hat\alpha^2-\hat\beta^2)^2}=\hat q^2$.
\end{proof}

\begin{remark}
Notice that in case of the SP equation we do not need to explore the compatibility of the differential equations \eqref{hat_hat_U} and \eqref{hat_hat_V} to obtain a representation of a solution of the SP equation in the $(y,t)$ variables \eqref{SP_in_y}.
\end{remark}

\begin{remark}
Notice that
\[
 \hat u=2(\hat\alpha\hat\beta_t-\hat\beta\hat\alpha_t).
\]
Thus $\hat u$, $\hat w$, $\hat q$, and $x(y,t)$ can be written in terms of $\hat\alpha$ and $\hat\beta$ only.
\end{remark}

\begin{remark}\label{rem:mean-zero-RH}
For a smooth periodic reconstruction, the zero mean follows from the RH coefficient identities. Indeed,
\[
 \hat u_{yt}
 =2(\hat\alpha_t\hat\beta+\hat\alpha\hat\beta_t)
 =\hat u(\hat\alpha^2-\hat\beta^2)=\frac{\hat u}{\hat q}.
\]
Consequently,
\[
 \int_0^\theta\frac{\hat u(y,t)}{\hat q(y,t)}\dd y
 =\hat u_t(\theta,t)-\hat u_t(0,t)=0.
\]
Since $\dd x=\dd y/\hat q$, this is the zero mean in the original spatial variable.
\end{remark}

Now let us discuss sufficient conditions on the contour $\Gamma$, jump matrix $J_0(k)$ and residue conditions that provide the symmetries \eqref{sym_1} and \eqref{sym_2}.

\begin{proposition}\label{prop:sym}
    Let the contour $\Gamma$ be invariant w.r.t. $k\mapsto-\bar k$ and $k\mapsto\bar k$, the jump matrix $J_0$ satisfy the symmetries
    \begin{equation}
 \label{sym_J_1}
 J_0(k)=\overline{J_0(-\bar k)}
 \end{equation}
and 
 \begin{equation}
 \label{sym_J_2}
J_0^{-1}(k)=\sigma_2 \overline{J_0(\bar k)}\sigma_2,
 \end{equation} 
 and the parameters of the residue conditions satisfy \eqref{res-sym}. Assume also that the first-column pole set is invariant under $k\mapsto-\bar k$ and
$c_{-\bar k_j}=-\overline{c_{k_j}}$, with the corresponding conditions for the second-column poles. The reflection $k\mapsto-\bar k$ is understood to preserve the plus and minus sides, while $k\mapsto\bar k$ interchanges them.

    Then the solution $\hat M$ of the RH problem  \eqref{jump_cond}--\eqref{res_hatM__loc} satisfies symmetries \eqref{sym_1} and \eqref{sym_2}
\end{proposition}


\section{RH problem formalism for the periodic SP problem in $y,t$ variables }\label{sec:6}

In Section \ref{sec:5} we discussed how a function solving, locally, the SP equation 
in $y,t$ variables appears from the solution of a Riemann-Hilbert
problem parametrized by $y$ and $t$ in an appropriate way. In this section, we 
present a particular  Riemann-Hilbert problem associated with the periodic problem
for the SP equation in $y,t$ variables in the sense (i) that  it can be formulated in terms
of the spectral functions determined by the initial data only, (ii) whose solution
evaluated at $t=0$ recovers the initial data, and (iii) whose solution provides
the periodicity, i.e.  ensures that the values of the solution 
of the SP equation at the ends of the interval match.

\subsection{Construction of the master RH problem} It is  the pre-RH problem associated with the problem on the interval that \emph{suggests}
considering the RH problem with the jump conditions (i) parametrized by $y$ and $t$
and (ii) given in terms of the spectral functions $a(k)$, $b(k)$, 
and \emph{ a function ${\mathfrak R}(k)$ that can be calculated in terms of $a(k)$ and $b(k)$}.

Namely,  in the case of the periodic problem, $\Phi_{\infty 4}(L,t,k)=\Phi_{\infty 1}(0,t,k)$ (they satisfy the same integral equation since $V_{\infty}(L,t,k)=V_{\infty }(0,t,k)$) which implies
$B_1(k)=B(k)$ and $A_1(k)=A(k)$.
It follows that in this case, the l.h.s. of the Global Relations 
\eqref{gl_rel} and \eqref{gl_rel_til}, after dividing respectively by
$A(k)A_1(k)$ and $\tilde A(k)\tilde A_1(k)$, becomes 
 quadratic expressions w.r.t. respectively $R(k):=\frac{B(k)}{A(k)}$ and $\tilde R(k):=\frac{\tilde B(k)}{\tilde A(k)}$,
 with coefficients expressed in terms of $a$, $b$, $\tilde a$, and $\tilde b$.
 The decisive step consists in replacing the r.h.s. in the Global relations by $0$,
 which gives two quadratic equations, respectively for $\mathfrak R$ and $\tilde{\mathfrak R}$:
\begin{equation}\label{calR}
    \mathfrak{R}^2(k)b^*( k)\eul^{2ik\theta}+\mathfrak{R}(k)(a^*( k)\eul^{2ik\theta}-a(k))+b(k)=0.
\end{equation}
        and 
\begin{equation}\label{caltilR}
{\tilde{\mathfrak R}}^2(k)\tilde b^*(k)e^{2ikL}
		+ {\tilde{\mathfrak R}}(k) \left(e^{2ikL} \tilde a^*(k) - \tilde a(k)\right)  +\tilde b(k) = 0.
\end{equation}

         Taking into consideration that $\tilde a(k)=a(k) e^{ik(L-\theta)}$ and 
         $\tilde b(k)=b(k) e^{ik(L-\theta)}$, the latter equation reduces to \eqref{calR},
         which will be taken as the (part of) \emph{definitions} of functions 
         $\mathfrak R(k)$ and ${\tilde{\mathfrak R}}(k)$.

Being solutions of a  quadratic equation, ${\mathfrak R}(k)$ and $\tilde{\mathfrak R}(k)$ 
 can be expressed analytically by 
\begin{equation}\label{R12}
\mathfrak{R}(k)=\frac{a(k)\eul^{-ik\theta}-a^*( k)\eul^{ik\theta}\pm\sqrt{\Delta^2(k)-4}}{2\eul^{ik\theta}b^*( k)}
\end{equation}
 and similarly for $\tilde{\mathfrak R}(k)$, where $\Delta(k)=a(k)\eul^{-ik\theta}+\overline{a(\bar k)}\eul^{ik\theta}$ and where  
		in order to fix a one-valued branches of ${\mathfrak R}(k)$  and 
        $\tilde{\mathfrak R}(k)$ (which can be different) for $k$ in the complex plane, one needs to introduce cuts
        connecting 
		simple zeros of $\Delta^2({z})-4$ and to specify the values 
        of ${\mathfrak R}(k)$  and 
        $\tilde{\mathfrak R}(k)$ at dedicated  points.      
        Looking at the r.h.s. of the Global Relations suggests to
specify them as follows:
\begin{equation}\label{spec_R}
{\mathfrak R}(\infty)=0, \quad \tilde{\mathfrak R}\left(0\right)=0.    
\end{equation}

Notice that there are two possibilities: either 
${\mathfrak R}\left(0\right)=0$ and thus
${\mathfrak R}(k)\equiv\tilde{\mathfrak R}(k)$ or 
${\mathfrak R}\left(0\right)\ne 0$ and thus
${\mathfrak R}(k)$ and $\tilde{\mathfrak R}(k)$  are the different branches of the same analytic function defined on the corresponding Riemann surface.

Before passing to the formulation of the RH problem with the jump given in terms of 
$a(k)$, $b(k)$, ${\mathfrak R}(k)$,  and $\tilde{\mathfrak R}(k)$, we present the analytic 
properties of ${\mathfrak R}(k)$  and $\tilde{\mathfrak R}(k)$ following their definitions.
But before doing that, we notice that the scattering matrix $s(k)$ is related to the monodromy matrix. Let $F(\xi,k)$ be the solution of the original spatial equation \eqref{Lax-x}, considered at $t=0$, fixed by $F(0,k)=I$. Since $P(L,0)=P(0,0)$, the monodromy of the gauged equation is
\[
 \mathcal M(k)=P(0,0)F(L,k)P^{-1}(0,0)
 =\eul^{-\ii k\theta\sigma_3}s^{-1}(k)
 =
 \begin{pmatrix}
 a(k)\eul^{-\ii k\theta}&-b(k)\eul^{-\ii k\theta}\\
 b^*(k)\eul^{\ii k\theta}&a^*(k)\eul^{\ii k\theta}
 \end{pmatrix}.
\]
In particular, $\Delta(k)=\operatorname{tr}\mathcal M(k)=\operatorname{tr}F(L,k)$.

\begin{proposition}

\begin{enumerate}[1.]
\item $\overline{\Delta(\bar k)}=\Delta(k)$ and $\overline{\Delta(-\bar k)}=\Delta(k)$. In particular, if $\xi_j$ is a zero of $\Delta^2(k)-4$, the symmetries imply that $\bar \xi_j$, $-\xi_j$ and $-\bar \xi_j$ are also zeros of $\Delta^2(k)-4$.

Combining these symmetries, we also get $\Delta(k)=\Delta(-k)$.

    \item For $x\in \mathbb{R}$, we have $|\Delta (x)|\leq2$ and $\Delta (x)\in\mathbb{R}$.

    \item Let $x_j$ be a real zero of $g(k):=\Delta^2(k)-4$. Then
\begin{enumerate}
    \item  its multiplicity is even and strictly greater than $1$;

    \item $b(x_j)=0$ (i.e. it is a Dirichlet eigenvalue).
\end{enumerate}

\item For $\ii x\in\ii\mathbb R$, we have $\Delta(\ii x)\in\mathbb R$. Moreover,
\begin{equation}\label{delta-imaginary-bound}
 \Delta(\ii x)\ge 2\cosh(xL)>2,\qquad x\ne0.
\end{equation}
In particular, $\Delta^2(k)-4$ has no nonzero purely imaginary zeros.

\item $\Delta^2(k)-4$ has a zero of even multiplicity (in fact, multiplicity two) at $k=0$.

\end{enumerate}

\end{proposition}

\begin{proof} \noindent 1. The first symmetry follows from the definition of $\Delta(k)$, while the second symmetry follows from the corresponding symmetry of $a(k)$ (see \eqref{symm_a}).

\noindent 2. First, notice that $\overline{\Delta(\bar k)}=\Delta(k)$ implies that for $x\in\mathbb{R}$ we have $\overline{\Delta(x)}=\Delta(x)$ and thus it is real. Moreover,
\begin{equation}
\label{del_re}\Delta(x)=a(x)\eul^{-\ii x\theta} + \overline{a(x)}\eul^{\ii x\theta} =2\Re\left(a(x)\eul^{-\ii x\theta}\right)
   \end{equation}
Thus, \begin{equation}\label{del_re_2} 
|\Delta(x)|^2\leq4|a(x)|^2.
\end{equation}
On the other hand, the determinant condition implies 
\begin{equation}\label{det_re}
 |a(x)|^2+|b(x)|^2=1, \quad  |a(x)|^2\leq 1.
\end{equation}
Combining \eqref{del_re_2} and \eqref{det_re}, we get $|\Delta(x)|^2\leq4$.

\medskip
\noindent 3.
Notice that $\Delta(x_j)\neq 0$ and the direct computation gives $g'(k)=2\Delta(k)\Delta'(k)$. Thus, we need to prove that $\Delta'(x_j)=0$.

Combining \eqref{del_re} and $\Delta^2(x_j)=4$, we get 
\begin{equation}
\label{del_x_j}
\Delta(x_j)=2\Re\left(a(x_j)\eul^{-\ii x_j\theta}\right)=\pm 2
\end{equation}
Rewriting \eqref{det_re} (at $x_j$) as 
\[
a(x_j)\eul^{-\ii x_j\theta} \overline{a(x_j)\eul^{-\ii x_j\theta}}+b(x_j)\overline{b(x_j)}=1
\]
and using 
\eqref{del_x_j}, we obtain
\begin{equation}
    \label{del_x_j_2}  \Im\left(a(x_j)\eul^{-\ii x_j\theta}\right)=0, \quad b(x_j)=0.
\end{equation}

Differentiating the determinant relation w.r.t. $k$ we get
\[
0=a'(k)a^*(k)+a(k){a^*}'(k)+b'(k)b^*(k)+b(k){b^*}'(k)
\]
Now, substituting $k=x_j$ and using \eqref{del_x_j_2}, we get
\[
0=a'(x_j)\overline{a(x_j)}+a(x_j)\overline{a'(x_j)}=a'(x_j)\eul^{-\ii x_j\theta}\overline{a(x_j)\eul^{-\ii x_j\theta}}+a(x_j)\eul^{-\ii x_j\theta}\overline{a'(x_j)\eul^{-\ii x_j\theta}}.
\]
Thus, using \eqref{del_x_j_2} again, we get
\begin{equation}
\label{der_del_x_j}
   0=a'(x_j)\eul^{-\ii x_j\theta}+\overline{a'(x_j)\eul^{-\ii x_j\theta}}=2\Re\left(a'(x_j)\eul^{-\ii x_j\theta}\right). 
\end{equation}

Differentiating $\Delta (k)$ w.r.t. $k$, we get
\[
\Delta '(k)=a'(k)\eul^{-\ii k \theta}+\overline{a'(\bar k)}\eul^{\ii k \theta}+\ii\theta\left(-a(k)\eul^{-\ii k\theta}+\overline{a(\bar k)}\eul^{\ii k\theta}\right)
\]

Finally, substituting $k=x_j$ we get
\[
 \Delta'(x_j)
 =2\Re\left(a'(x_j)\eul^{-\ii x_j\theta}\right)
 +2\theta\Im\left(a(x_j)\eul^{-\ii x_j\theta}\right).
\]
Thus, \eqref{del_x_j_2} and \eqref{der_del_x_j} imply that $\Delta'(x_j)=0$. Moreover, $g$ is real analytic and $g(x)\leq0$ for real $x$. Therefore every real zero has even multiplicity.

\medskip
\noindent 4. The symmetry $\overline{\Delta(-\bar k)}=\Delta(k)$ implies that for $\ii x\in\ii\mathbb R$ we have $\overline{\Delta(\ii x)}=\Delta(\ii x)$ and thus it is real. Moreover,
\[
 \Delta(\ii x)=a(\ii x)\eul^{x\theta}+a(-\ii x)\eul^{-x\theta}.
\]

To prove \eqref{delta-imaginary-bound}, write $k=\ii\eta$, $\eta\in\mathbb R\setminus\{0\}$, and use the original spatial Lax equation. With $w_0=u_{0x}$, it reads
\begin{equation}\label{imaginary-spatial-system}
 F_\xi(\xi,\ii\eta)=\eta A(\xi)F(\xi,\ii\eta),
 \qquad
 A(\xi)=\begin{pmatrix}1&w_0(\xi)\\w_0(\xi)&-1\end{pmatrix},
 \qquad F(0,\ii\eta)=I.
\end{equation}
Here $T(\ii\eta):=F(L,\ii\eta)$ is real, $\det T(\ii\eta)=1$, and $\Delta(\ii\eta)=\operatorname{tr}T(\ii\eta)$. Since
\[
 A^T\sigma_3+\sigma_3A=2I,
\]
every vector solution $\phi$ of \eqref{imaginary-spatial-system} satisfies
\begin{equation}\label{imaginary-coercive-identity}
 \frac{\dd}{\dd\xi}\left(\phi^\dagger\sigma_3\phi\right)
 =2\eta\|\phi\|^2,
\end{equation}
where $\dagger$ denotes the Hermitian transpose.

Let $\rho$ be an eigenvalue of $T(\ii\eta)$ and $v\ne0$ a corresponding eigenvector. Applying \eqref{imaginary-coercive-identity} to $\phi(\xi)=F(\xi,\ii\eta)v$ gives
\begin{equation}\label{imaginary-floquet-identity}
 (|\rho|^2-1)v^\dagger\sigma_3v
 =2\eta\int_0^L\|\phi(\xi)\|^2\dd\xi.
\end{equation}
The right-hand side is nonzero. Hence no Floquet multiplier lies on the unit circle. Since $T(\ii\eta)$ is real and has determinant one, this implies $|\Delta(\ii\eta)|>2$. Continuity and $\Delta(0)=2$ then imply $\Delta(\ii\eta)>2$ on both half-axes.

For the quantitative bound, let $\eta>0$ and choose the expanding multiplier $\rho>1$ and a real eigenvector $v$. Equation \eqref{imaginary-floquet-identity} gives $h(0)>0$, where $h(\xi)=\phi(\xi)^T\sigma_3\phi(\xi)$. Moreover,
\[
 h'(\xi)=2\eta\|\phi(\xi)\|^2\ge2\eta h(\xi).
\]
Consequently,
\[
 \rho^2h(0)=h(L)\ge\eul^{2\eta L}h(0),
 \qquad \rho\ge\eul^{\eta L}.
\]
It follows that
\[
 \Delta(\ii\eta)=\rho+\rho^{-1}\ge2\cosh(\eta L)>2.
\]
The evenness of $\Delta$ gives the same bound for $\eta<0$.

\medskip
\noindent 5. First notice that $\Delta(k)=\Delta(-k)$ implies that $\Delta(k)$ (and, thus, $\Delta^2(k)-4$) is even. This implies that if $k=0$ is zero of $\Delta^2(k)-4$, it has even multiplicity.

Combining \eqref{s_i2}
and \eqref{s_til_s}, we get
\begin{equation}
\label{a_at_0}
    a(k)=1+\ii a_0k+O(k^2), \quad
    a^*(k)=1-\ii a_0k+O(k^2),
\end{equation}
where $a_0=\theta-L$, and thus $\Delta(0)=2$. More precisely, expanding the original spatial equation gives
\[
 \Delta(k)=2-k^2\left(L^2+\left(\int_0^L w_0(\xi)\dd\xi\right)^2\right)+O(k^4).
\]
Since $u_0$ is periodic, $\int_0^Lw_0(\xi)\dd\xi=0$. Therefore
\begin{equation}
 \label{delta_at_0}
 \Delta^2(k)-4=-4L^2k^2+O(k^4).
\end{equation}

\end{proof}

We will additionally assume that all nonreal zeros of $\Delta^2(k)-4$ (main eigenvalues) are simple and that there is a finite number of such zeros. The absence of nonzero purely imaginary zeros and the even multiplicity of all real zeros have already been proved above.

For the RH formulations below, we also impose the usual nondegeneracy assumptions: poles of $\mathfrak R$ in the exterior region are simple, finite in number, and disjoint from the jump contour and branch endpoints; $\tilde{\mathfrak R}$ is holomorphic in the small disk; the small circle avoids all branch points and spectral poles; and the denominators used at branch endpoints do not vanish. A possible pole of $\mathfrak R$ at $0$ is harmless, since $\mathfrak R$ is used only outside the small disk. The zeros of $a$ used in the initial-data RH problem are assumed simple, off the real axis, and disjoint from the branch endpoints and ratio poles. Coincident or contour singularities require a separate limiting formulation.

\begin{proposition}\label{prop:R}

(1) 
Let $\mathfrak{R}_{(1)}$ and $\mathfrak{R}_{(2)}$ be two solutions of \eqref{calR}. Then the following properties hold:

  \begin{equation}\label{propR1}
        \mathfrak{R}_{(1)}\mathfrak{R}_{(2)}=-\frac{\mathcal M_{12}}{\mathcal M_{21}}=\frac{b(k)\eul^{-2\ii k \theta}}{b^*( k)}=\frac{\tilde b(k)\eul^{-2\ii k L}}{{\tilde b^*(k)}}
    \end{equation}
        \begin{equation}\label{propR2}
        \mathfrak{R}_{(1)}+\mathfrak{R}_{(2)}=\frac{\mathcal M_{11}-\mathcal M_{22}}{\mathcal M_{21}}=\frac{a(k)\eul^{-\ii k \theta}-a^*( k)\eul^{\ii k \theta}}{b^*( k)\eul^{\ii k \theta}}=\frac{\tilde a(k)\eul^{-\ii k L}-{\tilde a^*(k)}\eul^{\ii k L}}{{\tilde b^*(k)}\eul^{\ii k L}}
    \end{equation}
 
            \begin{equation}\label{propR3}
        \mathfrak{R}_{(2)}^*=-\frac{1}{\mathfrak{R}_{(1)}}, \quad \mathfrak{R}_{(1)}^*=-\frac{1}{\mathfrak{R}_{(2)}}
    \end{equation}

    \begin{equation}\label{propR4}
        \mathfrak{R}_{(2)}\mathcal M_{21}=\mathfrak{R}_{(2)}^*\mathcal M_{12}, \quad \mathfrak{R}_{(1)}\mathcal M_{21}=\mathfrak{R}_{(1)}^*\mathcal M_{12}, \quad \mathfrak{R} \mathcal M_{21}=\mathfrak{R}^*\mathcal M_{12}
    \end{equation}

(2)  Let $\mathfrak{R}(k)$ be either $\mathfrak{R}_{(1)}(k)$ or $\mathfrak{R}_{(2)}(k)$.
Then
\begin{equation}\label{propR5}
      \mathfrak{R}(k) b^*( k)\eul^{\ii k \theta}=-\mathfrak{R}^*(k)b(k)\eul^{-\ii k \theta}
    \end{equation}
\begin{equation}\label{propR6}
        (a+b\mathfrak{R}^*)(a^*+b^*\mathfrak{R})=1
    \end{equation}
and thus
 \begin{equation}\label{propR7}
        (a+b\mathfrak{R}^*)\neq 0 \ \ \text{where both factors are finite}.
    \end{equation} 

(3) Let $\mathfrak{R}_\pm (k)$ be the limiting values of $\mathfrak{R}$ as its argument approaches 
           the (oriented) branch cut from the corresponding side. Then
\begin{equation}\label{propR8}
        \mathfrak{R}_+(k) \mathfrak{R}^{*}_-(k)=\mathfrak{R}_-(k) \mathfrak{R}^{*}_+(k)=-1
    \end{equation}  

        \begin{equation}\label{propR9}
        (a^*+b^*\mathfrak{R}_-)(k)=\eul^{-2\ii k \theta} (a+b\mathfrak{R}^{*}_+)(k),\quad (a+b\mathfrak{R}^{*}_-)(k)=\eul^{2\ii k \theta}(a^*+b^*\mathfrak{R}_+)(k)
    \end{equation}

\end{proposition}
Additionally, 
\begin{equation}\label{propR10}
        \mathfrak{R}(-k)=\overline{\mathfrak{R}(\bar k)}
    \end{equation}   
    due to the corresponding symmetries $a(k)$ and $b(k)$.

The RH problem formalism for the periodic problem for the SP equation involves the following
steps.

\noindent\textbf{Step 1}. Formulate the RH problem for $\hat M(y,t,k)$, depending on $y$ and $t$ as parameters,
providing
\begin{itemize}
    \item Jump conditions
    \begin{equation}\label{hatM-jump}
        \hat M_+(y,t,k)= \hat M_-(y,t,k)\hat J(y,t,k), \quad k\in \hat\Sigma
    \end{equation}
    replacing in \eqref{jump-xt} and \eqref{j0-xt}
    $p(x,t,k)$ by $\hat p(y,t,k)=y-\frac{t}{4k^2}$,
     $R_1(k)$ by $\mathfrak R(k)$, and $\tilde R_1(k)$ 
    by $\tilde{\mathfrak R}(k)$;
    
    \item
    Normalization condition $\hat M(k)\to I$ as $k\to\infty$, together with boundedness at $0$ and the finite contour vertices as in Section \ref{sec:5};

\item Residue conditions at poles $\mu_j$ of $\mathfrak R(k)$ in $\mathbb{C}^+\setminus \{|k|<\epsilon\}$ 
and at poles $\bar\mu_j$ of $\mathfrak R^*(k)$ in $\mathbb{C}^-\setminus \{|k|<\epsilon\}$:
\begin{subequations}\label{res_hatM_R}
    \begin{align}
    \label{res_hatM-a}
\Res_{\mu_j}\hat M^{(2)}&= \hat M^{(1)}(\mu_j) \eul^{-2\ii \mu_j \hat p(\mu_j)} \Res_{\mu_j} \mathfrak R, ~ \mu_j\in \mathbb{C}^+\setminus \{|k|<\epsilon\},\\
\label{res_hatM-b}
\Res_{\bar\mu_j}\hat M^{(1)}&=-\hat M^{(2)}(\bar\mu_j) \eul^{2\ii \bar\mu_j \hat p(\bar\mu_j)}\Res_{\bar\mu_j} {\mathfrak R}^*,~ \bar\mu_j\in \mathbb{C}^-\setminus \{|k|<\epsilon\}. 
\end{align}
\end{subequations}

\end{itemize}

\noindent\textbf{Step 2}. Prove that the solution of the RH problem above,
being evaluated at $t=0$, gives rise to the initial condition $\hat u(y,0)$.

\smallskip

\noindent\textbf{Step 3}. Prove that the solution of the RH problem above
gives rise to the solution of the SP equation periodic in $y$ for all fixed $t>0$.

Now let's make the jump conditions \eqref{hatM-jump} precise.
First, we notice that since we are going to replace analytic functions (in the pre-RH problem) by ones having jumps across the branch cuts (${\mathfrak R}(k)$ and 
${\mathfrak R}^*(k)$), the resulting jump contour $\hat\Sigma$ has to include the branch cuts:
$\hat\Sigma =\mathbb{R}\cup\{|k|=\epsilon\}\cup_{j=1}^{N}\Sigma_j$, where and
$\{\Sigma_j\}$ are segments connecting 
simple zeros of $\Delta^2({k})-4$.

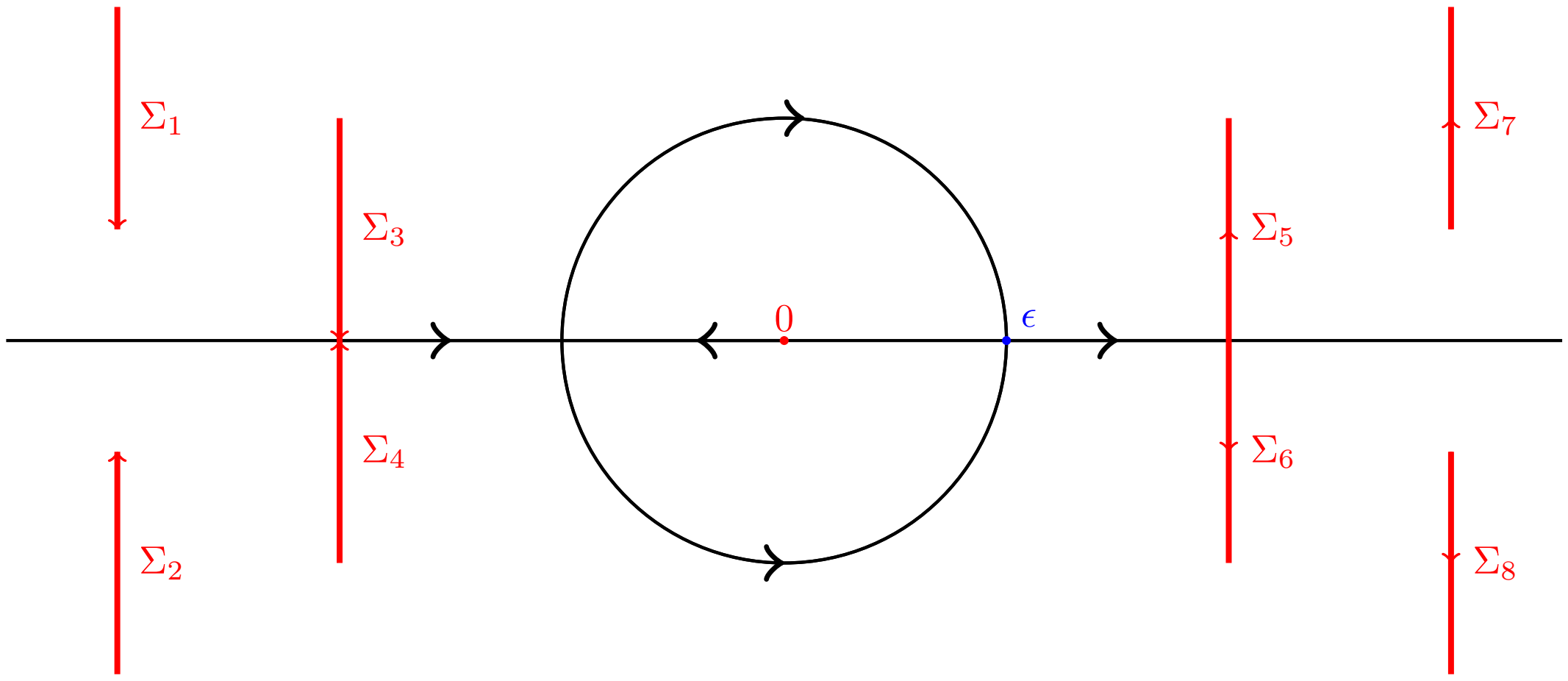
\begin{figure}[ht]
    \centering
\begin{tikzpicture}[scale=2]
\draw[red] (0,0) circle (1);
\draw[thick] (-3.5,0) -- (3.5,0);
\draw[thick, postaction={decorate},
      decoration={markings, mark=at position 0.5 with {\arrow[scale=2]{>}}}]
      (-2,0) -- (-1,0);
\draw[thick, postaction={decorate},
      decoration={markings, mark=at position 0.7 with {\arrow[scale=2]{<}}}]
      (-1,0) -- (0,0);  
\draw[thick, postaction={decorate},
      decoration={markings, mark=at position 0.5 with {\arrow[scale=2]{>}}}]
      (1,0) -- (2,0);  
\draw[thick, postaction={decorate},
      decoration={markings, mark=at position 0.25 with {\arrow[scale=2]{<}}}] 
      (0,0) circle (1);
\draw[thick, postaction={decorate},
      decoration={markings, mark=at position 0.75 with {\arrow[scale=2]{>}}}] 
      (0,0) circle (1);
\node at (0,0.1) {\textcolor{red}{$0$}};
\fill[red] (0,0) circle(0.02);
\node at (1.1,0.1) {\textcolor{blue}{$\epsilon$}};
\fill[blue] (1,0) circle(0.02);

\draw[line width=1.5pt, red] (-2, 0) -- (-2, -1);
             \draw[thick, red, <-] (-2, 0) -- (-2, -0.5);
\draw[line width=1.5pt, red] (-2, 0) -- (-2, 1);
             \draw[thick, red, <-] (-2, 0) -- (-2, 0.5);

\draw[line width=1.5pt, red] (2, 0) -- (2, -1);
             \draw[thick, red, ->] (2, 0) -- (2, -0.5);
\draw[line width=1.5pt, red] (2, 0) -- (2, 1);
             \draw[thick, red, ->] (2, 0) -- (2, 0.5);

\draw[line width=1.5pt, red] (-3, -0.5) -- (-3, -1.5);
             \draw[thick, red, <-] (-3, -0.5) -- (-3, -1);
\draw[line width=1.5pt, red] (-3, 0.5) -- (-3, 1.5);
             \draw[thick, red, <-] (-3, 0.5) -- (-3, 1);

\draw[line width=1.5pt, red] (3, -0.5) -- (3, -1.5);
             \draw[thick, red, ->] (3, -0.5) -- (3, -1);
\draw[line width=1.5pt, red] (3, 0.5) -- (3, 1.5);
             \draw[thick, red, ->] (3, 0.5) -- (3, 1);

 \node at (-2.8, 1) {$\red{\Sigma_1}$};      \node at (-2.8, -1) {$\red{\Sigma_2}$}; 

  \node at (-1.8, 0.5) {$\red{\Sigma_3}$};      \node at (-1.8, -0.5) {$\red{\Sigma_4}$}; 

   \node at (2.2, 0.5) {$\red{\Sigma_5}$};      \node at (2.2, -0.5) {$\red{\Sigma_6}$};   

 \node at (3.2, 1) {$\red{\Sigma_7}$};      \node at (3.2, -1) {$\red{\Sigma_8}$};   

\end{tikzpicture}
    \caption{The oriented contour $\hat\Sigma$}
    \label{fig:hatSigma}
\end{figure}

A typical $\hat\Sigma$ is shown in Figure \ref{fig:hatSigma}, where the small circle (of radius 
$\epsilon$) is chosen such that $\{\Sigma_j\}$
do not intersect the corresponding disk.

Then $\hat J(y,t,k) = \ee^{-ik\hat p \sigma_3}\hat J_0(k)\ee^{ik\hat p \sigma_3}$,
where $\hat J_0(k)$ is as follows.

\noindent\textbf{(i) $k\in\mathbb{R}\cup\{|k|=\epsilon\}$}. Here the construction of $\hat J_0(k)$
follows that of $J_0(k)$ \eqref{j0-xt}, 
 with 
$\Gamma_1(k)$ and $\tilde\Gamma_1(k)$ replaced by respectively  $G_1(k)$ and
$\tilde G_1(k)$,  and $\Gamma_2(k)$ replaced by   $G_2(k)$:
\begin{equation}\label{j0-yt}
    \hat J_0(k)=
    \begin{cases}
\begin{pmatrix}
    1& G_2(k)\\
   0 &1
\end{pmatrix}, & \Im k>0, |k|=\epsilon,\\

\begin{pmatrix}
          1& 0\\
         \frac{G^*_1(k)}{a^{*2}(k)}&1  
\end{pmatrix}
\begin{pmatrix}
          1& r(k)\\
         r^*(k)&1  + r(k) r^*(k)
\end{pmatrix}
\begin{pmatrix}
          1& \frac{G_1(k)}{a^{2}(k)}\\
         0&1  
\end{pmatrix}
, &  \Im k = 0, |k|>\epsilon,\\
\begin{pmatrix}
          1& -\frac{\tilde G_1(k)}{\tilde a^{2}(k)}\\
         0&1  
\end{pmatrix}
\begin{pmatrix}
          1 + \tilde r(k) \tilde r^*(k)& -\tilde r(k)\\
         -\tilde r^*(k)&1  
\end{pmatrix}
\begin{pmatrix}
          1& 0\\
         -\frac{\tilde G^*_1(k)}{\tilde a^{*2}(k)}&1   
\end{pmatrix}, &  \Im k = 0, |k|<\epsilon,\\
\begin{pmatrix}
    1& 0\\
   G_2^*(k) &1
\end{pmatrix}, & \Im k<0, |k|=\epsilon.

\end{cases}
\end{equation}

Here
\begin{align}\label{G-G1}
    G_1&:=\frac{\ee^{2\ii k \theta}a{\mathfrak R}}{a+b{\mathfrak R}^* }=a^2({\mathfrak R}-\frac{b}{a}),\quad
     \tilde G_1:=\frac{\ee^{2\ii k L} \tilde a \tilde{\mathfrak R}}{ \tilde a+ \tilde b \tilde{\mathfrak R}^* }= \tilde a^2( \tilde{\mathfrak R}-\frac{ \tilde b}{ \tilde a}),\\
     \label{G2}
      G_2&:=\ee^{\ii k (L+\theta)}\frac{{\mathfrak R}-\tilde{\mathfrak R}}{(a-b^*{\mathfrak R}\ee^{2\ii k \theta})(\tilde a-\tilde b^*\tilde{\mathfrak R}\ee^{2\ii k L})}=\begin{cases}
          0, & \text{if}\ \mathfrak R = \tilde{\mathfrak R},\\
          ({\mathfrak R}-\tilde{\mathfrak R}), & \text{if}\ \mathfrak R \ne \tilde{\mathfrak R}.
      \end{cases}
\end{align}
(we have used properties of ${\mathfrak R}$ and $\tilde{\mathfrak R}$
from Proposition \ref{prop:R}).

\medskip

\noindent\textbf{(ii)}
On $\{\Sigma_j\}$, the jumps are given in terms of the jumps of ${\mathfrak R}$; namely,
	\begin{equation}\label{sigma_j_ver}
	    \hat J_0(k)=\begin{cases} \begin{pmatrix}
		1  & {\mathfrak R}_+(k)-{\mathfrak R}_-(k) \\ 0 & 1
	\end{pmatrix}, & k \in \Sigma_j\cap \D{C}^+\\
	\begin{pmatrix}
		1  & 0 \\
		{{\mathfrak R}^*}_-(k)-{{\mathfrak R}^*}_+(k) & 1
	\end{pmatrix}, & k \in \Sigma_j\cap \D{C}^-
	\end{cases}
	\end{equation}

\begin{remark}\label{rem:master-solvability}
Notice that
\[
 (\mathfrak R_+-\mathfrak R_-)^*
 =\mathfrak R^*_- -\mathfrak R^*_+.
\]
Together with the contour orientations, this gives the required Schwarz reflection condition on the nonreal parts of the contour. On the real axis, the identities \eqref{G-G1} simplify the jump to
\[
 \hat J_0(k)=
 \begin{cases}
 \begin{pmatrix}1&\mathfrak R\\\mathfrak R^*&1+\mathfrak R\mathfrak R^*\end{pmatrix},
     &|k|>\epsilon,\\[2mm]
 \begin{pmatrix}1+\tilde{\mathfrak R}\tilde{\mathfrak R}^*&-\tilde{\mathfrak R}\\
 -\tilde{\mathfrak R}^*&1\end{pmatrix},
     &|k|<\epsilon.
 \end{cases}
\]
Both matrices are Hermitian positive definite: their determinant is one and their upper-left entry is positive. The phase conjugation is unitary for real $k$, so $\hat J+\hat J^\dagger$ is positive definite there.

Under the stated nondegeneracy and local regularity assumptions, the jumps belong to $I+(L^2\cap L^\infty)$ on a finite piecewise smooth contour, with the indicated decay on its unbounded parts. The simple poles can be removed by triangular factors on small circles paired under the two symmetries; the conjugate residue conditions give the same Schwarz reflection condition on these circles. Zhou's vanishing lemma \cite{Z89}, together with the index-zero Fredholm formulation, then gives unique solvability in this RH class. 
\end{remark}

\subsection{Verification of initial conditions}
The idea of verification of initial conditions (Step 2) is as follows: transform the master RH problem
evaluated at $t=0$ to a RH problem  that characterizes the initial condition $\hat u_0(y)$
(or $\hat w_0(y)$) in such a way that does not affect the procedure of ``extracting''
$\hat w$ from the solution of the RH problem.

We begin with the formulation of the RH problem associated with the Lax equation 
\eqref{Lax_y_y} at $t=0$, i.e., with the coefficients expressed in terms of $\hat w_0(y)$.
This RHP can be derived from the direct problem on the whole line with the potential supported on $[0,L]$ (see \cite{BSL17}). 

\begin{proposition}\label{prop:M-y}
    Let $a(k)$, $b(k)$ ($k\in\mathbb R)$, and $\{k_j\}_1^N$ ($k_j\in\mathbb{C}^+$ -- zeros of $a(k)$) be the spectral
    data determined by $u_0(x)$, $x\in[0,L]$
    through \eqref{rel_inf_1} and \eqref{s}, where \eqref{rel_inf_1} is considered at $t=0$.
   Let $w_0(x)=u_{0x}(x)$ and $\hat w_0(y)=w_0(x(y))$, where 
   $y(x)=\int_0^x\sqrt{1+w_0^2(\xi)}\dd\xi$. Then $\hat w_0(y)$ 
   can be given in terms of the solution of the following RH problem: for all 
   $y\in [0,\theta]$, 
   find a piece-wise (w.r.t. $\mathbb R$) meromorphic, $2\times 2$-valued 
   function $M^{(y)}(y,k)$ such that
   \begin{enumerate}
      \item The limiting values of $M^{(y)}(y,k)$ as $k$ approaches $\mathbb R$ from the upper 
    and lower half-planes ($M_+$ and $M_-$ respectively) are related by the jump 
    condition
    \begin{equation}\label{M-y-jump}
        M^{(y)}_+(y,k) = M^{(y)}_-(y,k)J^{(y)}_{\mathbb R}(y,k)\ \ \text{with}\ 
       J^{(y)}_{\mathbb R}(y,k)= e^{-iky\sigma_3}J_{0 \mathbb R}(k)  
       e^{iky\sigma_3},\ \ 
               k\in \mathbb R,
    \end{equation}
    where $J_{0 \mathbb R}(k) $ is given by 
    \begin{equation}
        \label{j0}
    J_{0 \mathbb R}(k)=\begin{pmatrix}
       1+\check r(k)\check r^*(k)&\check r^*(k)\\\check r(k)&1
   \end{pmatrix}
        \end{equation}
    with $\check r(k)=\frac{b^*(k)}{a(k)}$;
    
\item $M^{(y)}(y,k)\to I$ as $k\to\infty$;
\item 
$M^{(y)(1)}$ has simple poles at 
$k_j$, $j=1,\dots,N$,
$M^{(y)(2)}$ has simple poles at $\bar k_j$, and the residue conditions hold:
\begin{equation}\label{M-y-res}
\begin{aligned}
    \Res_{k_j} M^{(y)(1)}(y,k) & = \frac{e^{2\ii k_j y}}{\dot a(k_j)b(k_j)}  M^{(y)(2)}(y,k_j),\\
    \Res_{\bar k_j} M^{(y)(2)}(y,k) & =- \frac{e^{-2\ii\bar k_j y}}{\dot a^*(\bar k_j)b^*(\bar k_j)} M^{(y)(1)}(y,\bar k_j);
\end{aligned}
\end{equation}
       \end{enumerate}
Evaluating the solution $ M^{(y)}(y,k)$ of this RH problem  at $ k=0$       
\begin{equation}\label{y-struc}
    M^{(y)}\left(y,0\right)=\begin{pmatrix}
           \hat \alpha_0(y)&\hat \beta_0(y)\\
           -\hat \beta_0(y)&\hat \alpha_0(y)
       \end{pmatrix},
	\end{equation}
the following representation of $\hat w_0(y)$ holds    
\begin{equation}\label{m_0_via_RH_x_par}
    \hat{w}_0(y)=\frac{2\hat \alpha_0(y)\hat \beta_0(y)}{\hat \alpha_0^2(y)-\hat \beta_0^2(y)}.
    \end{equation}
\end{proposition}

\begin{proposition}\label{prop:con_pr}
    The solutions $\hat M(y,0,k)$ and $M^{(y)}(y,k)$ of the RH problems 
    \eqref{hatM-jump}-\eqref{res_hatM_R} and \eqref{M-y-jump}-\eqref{M-y-res} respectively
    can be related as follows:
    \begin{equation}
        \label{y-ini}
         M^{(y)}(y,k)=\hat{M} (y,0,k)\begin{cases}
\begin{pmatrix}
    \frac{1}{a(k)}&0\\0&a(k)
\end{pmatrix}
\begin{pmatrix}
   1&-G_1(k)\eul^{-2iky}\\
 0&1  
\end{pmatrix},\quad 
k\in\{ \Im k > 0, |k|>\epsilon \},\\

\begin{pmatrix}
    \frac{1}{\tilde a(k)}&0\\0&\tilde a(k)
\end{pmatrix}\begin{pmatrix}
   1&-\tilde G_1(k)\eul^{-2iky}\\
 0&1  
\end{pmatrix}\begin{pmatrix}
    \eul^{\ii k(L-\theta)}&0\\0& \eul^{-\ii k(L-\theta)}
\end{pmatrix},
\quad k\in\{  \Im k > 0, |k|<\epsilon \},\\

\begin{pmatrix}
  a^*(k)&0\\0&\frac{1}{a^*(k)}
\end{pmatrix}
\begin{pmatrix}
   1&0\\
 G_1^*(k)\eul^{2iky}&1  
\end{pmatrix}
,\quad k\in\{  \Im k < 0, |k|>\epsilon\},\\

\begin{pmatrix}
  \tilde a^*(k)  &0\\0&\frac{1}{\tilde a^*(k)}
\end{pmatrix}
\begin{pmatrix}
   1&0\\
 \tilde G_1^*(k)\eul^{2iky}&1  
\end{pmatrix}\begin{pmatrix}
    \eul^{\ii k(L-\theta)}&0\\0 &\eul^{-\ii k(L-\theta)}
\end{pmatrix}
,\quad k\in\{  \Im k < 0, |k|<\epsilon\}.
\end{cases}
\end{equation}
Moreover, $\hat\alpha(y,0)=\hat\alpha_0(y)$ and $\hat\beta(y,0)=\hat\beta_0(y)$, and thus $\hat w(y,0)=\hat w_0(y)$.
\end{proposition}

\begin{proof}
    The proof consists in the demonstration that the r.h.s. of \eqref{y-ini}
    meets all requirements of the RH problem for  $M^{(y)}(y,k)$ 
    in Proposition \ref{prop:M-y}.

    \medskip

\noindent (1) The jump conditions match by construction.

\medskip

\noindent (2) Using
\[
 G_1(k)\eul^{-2\ii ky}
 =\eul^{2\ii k(\theta-y)}
   \frac{a(k)\mathfrak R(k)}{a(k)+b(k)\mathfrak R^*(k)},
\]
the choice $\mathfrak R(k)=O(k^{-1})$ and the large-$k$ spectral asymptotics show that the upper-half-plane triangular factor is $I+O(k^{-1})$ for $0\le y\le\theta$. Here $\eul^{2\ii k(\theta-y)}$ is bounded in $\mathbb C^+$. The lower-half-plane case follows by symmetry. Thus the normalization conditions match.

\medskip

\noindent (3)
Let's check that the residue conditions match. 

For $k\in \mathbb{C}^+\setminus\{|k|<\epsilon\}$ we have
that 
\[
\hat M (y,0,k)
\begin{pmatrix}
   \frac{1}{a(k)}&-\frac{G_1(k)}{a(k)}\eul^{-2iky}\\
 0&a(k)  
\end{pmatrix}=\left(\frac{1}{a(k)}\hat M^{(1)} (y,0,k), -\frac{G_1(k)}{a(k)}\eul^{-2iky}\hat M^{(1)} (y,0,k)+a(k)\hat M^{(2)} (y,0,k)     \right),\] 
in particular, taking into account $\Res_{\mu_j}G_1 = a^2(\mu_j)\Res_{\mu_j} \mathfrak R $ (see \eqref{G-G1}) the singularities of the second column at $k=\mu_j$ cancel.

On the other hand, at $k=k_j$, the first column is singular due to the singularity of $\frac{1}{a(k)}$, and, using $\frac{G_1(k_j)}{a(k_j)}=-b(k_j)$, the corresponding residue condition is
\[
 \Res_{k_j}M^{(y)(1)}
 =\frac{\hat M^{(1)}(y,0,k_j)}{\dot a(k_j)}
 =\frac{\eul^{2\ii k_jy}}{\dot a(k_j)b(k_j)}
 M^{(y)(2)}(y,k_j),
\]
which is \eqref{M-y-res}.  

The singularities in $\mathbb{C}^-$  can be treated in a similar way.

Thus, by uniqueness of solution of RH problem \eqref{M-y-jump}--\eqref{M-y-res}, we conclude that \eqref{y-ini} holds.

\medskip

Finally, notice that $\tilde b\left(0\right) 
= \tilde{\mathfrak R}\left(0\right) =0$ in $\mathbb{C}^+$, which implies that $\tilde G_1\left(0\right)=0$.
Therefore, $\hat\alpha(y,0)=\hat\alpha_0(y)$ and $\hat\beta(y,0)=\hat\beta_0(y)$.

\end{proof}

\begin{remark}\label{rem:const_t_0}
    Notice that \eqref{y-ini} allows us to explicitly construct the solution $\hat M(y,t,k)$ of RH problem \eqref{hatM-jump}-\eqref{res_hatM_R} at $t=0$ from the solution $M^{(y)}(y,k)$ of RH problem \eqref{M-y-jump}-\eqref{M-y-res}, whose existence is guaranteed by construction. 
\end{remark}

\subsection{Verification of periodicity}
The idea of verification of periodicity (Step 3) is as follows: transform 
two RH problems, parametrized by $t$, which are 
the master RH problem for respectively $y=0$ and for $y=\theta$,
to two RH problems having the same data.

For this purpose we introduce

\begin{equation*}
    \hat P(t,k)=\begin{cases}
\begin{pmatrix}
   a^*+b^*\mathfrak R&0\\
 -b^*\eul^{2\ii k\theta}\eul^{-\frac{\ii  t}{2 k}}&a+b\mathfrak R^*  
\end{pmatrix}, & 
 k\in\{ \Im k > 0, |k|>\epsilon\},\\
\begin{pmatrix}
   a^*+b^*\tilde {\mathfrak R}&0\\
 -b^*\eul^{2\ii k\theta}\eul^{-\frac{\ii  t}{2 k}}&a+b\tilde {\mathfrak R}^*  
\end{pmatrix}, & 
 k\in\{ \Im k > 0, |k|<\epsilon\},\\
\begin{pmatrix}
   a^*+b^*\tilde {\mathfrak R}& b\eul^{-2\ii k\theta}\eul^{\frac{\ii  t}{2 k}}\\
0&a+b\tilde {\mathfrak R}^*  
\end{pmatrix}, & 
k\in\{ \Im k < 0, |k| <\epsilon\},\\
\begin{pmatrix}
   a^*+b^*\mathfrak R& b\eul^{-2\ii k\theta}\eul^{\frac{\ii  t}{2 k}}\\
0&a+b\mathfrak R^*  
\end{pmatrix}, & 
k\in\{ \Im k < 0, |k|>\epsilon\}
\end{cases}
\end{equation*}
and define $M^{(t)}(t,k):=\hat{M}(0,t,k)$ and 
$\hat M^{(t)}(t,k):=\hat{M}(\theta,t,k)\hat P(t,k)$. 
Our goal is to show that $M^{(t)}$ and $\hat M^{(t)}$ satisfy the same
conditions characterizing the respective RH problem.

\begin{proposition}\label{prop:jump}
Let $\Sigma=\hat\Sigma$ be oriented as in Figure \ref{fig:hatSigma}, and define the jump matrices $J^{(t)}(t,k)$ and $\hat J^{(t)}(t,k)$ by the relations
\[
M_+^{(t)}(t,k) = M_-^{(t)}(t,k)J^{(t)}(t,k), \quad 
\hat M_+^{(t)}(t,k) =\hat M_-^{(t)}(t,k)\hat J^{(t)}(t,k),
\]
for $k \in \Sigma$. Then the two jump matrices coincide. Moreover, \[J^{(t)}(t,k) = \hat J^{(t)}(t,k)=\ee^{\frac{\ii t}{4k}\sigma_3} \hat J_0(k) 
\ee^{-\frac{\ii t}{4k}\sigma_3}\] with
$\hat J_0(k)$ defined in \eqref{j0-yt}--\eqref{sigma_j_ver}.
\end{proposition}

\begin{proof}

By the definitions of $M^{(t)}$ and $\hat M^{(t)}$, we have
\begin{align}
    M^{(t)}_+&=\hat M_+(0)=\hat M_-(0)\hat J(0)= M^{(t)}_-\hat J(0),\\
    \hat M^{(t)}_+&=\hat M_+(\theta)\hat P_+=\hat M_-(\theta)\hat J(\theta)\hat P_+= \hat M^{(t)}_-\hat P_-^{-1}\hat J(\theta)\hat P_+,
\end{align}
and, thus, we need to check that 
\begin{equation}\label{per_jump}
   \hat P_-^{-1}(t,k)\hat J(\theta,t,k)\hat P_+(t,k)=\hat J(0,t,k),
\end{equation}
where $\hat J(0,t,k)=\ee^{\frac{\ii t}{4k}\sigma_3}\hat J_0(k) 
\ee^{-\frac{\ii t}{4k}\sigma_3}$ and 
$\hat J(\theta,t,k) = \ee^{-ik (\theta-\frac{t}{4k^2})\sigma_3 }\hat J_0(k) 
\ee^{ik (\theta-\frac{t}{4k^2})\sigma_3 }$ with $\hat J_0(k)$
given by \eqref{j0-yt} and \eqref{sigma_j_ver}.

\noindent (1) For $k\in  \{ \Im k > 0, |k|=\epsilon\}$, we have
\begin{align}
\hat P_-^{-1}&=        
\begin{pmatrix}
  a+b\tilde {\mathfrak R}^* &0\\
 b^*\eul^{2\ii k\theta}\eul^{-\frac{\ii  t}{2 k}}&a^*+b^*\tilde {\mathfrak R}   
\end{pmatrix},\quad \hat P_+= \begin{pmatrix}
   a^*+b^*\mathfrak R&0\\
 -b^*\eul^{2\ii k\theta}\eul^{-\frac{\ii  t}{2 k}}&a+b\mathfrak R^*  
\end{pmatrix}.
\end{align}
Then, using \eqref{det_rel}, \eqref{calR}, \eqref{propR5} and \eqref{propR6}, straightforward calculations show \eqref{per_jump}.

\noindent (2) For $k\in  \{ \Im k = 0, |k|>\epsilon\}$, we have
\begin{align}
\hat P_-^{-1}&=        
\begin{pmatrix}
 a+b\mathfrak R^*   & -b\eul^{-2\ii k\theta}\eul^{\frac{\ii  t}{2 k}}\\
0&a^*+b^*\mathfrak R 
\end{pmatrix},\quad \hat P_+= \begin{pmatrix}
   a^*+b^*\mathfrak R&0\\
 -b^*\eul^{2\ii k\theta}\eul^{-\frac{\ii  t}{2 k}}&a+b\mathfrak R^*  
\end{pmatrix}.
\end{align}
Then, using \eqref{calR} and \eqref{propR5}, straightforward calculations show \eqref{per_jump}.

\noindent (3) For $k\in  \{ \Im k = 0, |k|<\epsilon\}$, we have
\begin{align}
\hat P_-^{-1}&=        
\begin{pmatrix}
  a+b\tilde {\mathfrak R}^* &0\\
 b^*\eul^{2\ii k\theta}\eul^{-\frac{\ii  t}{2 k}}&a^*+b^*\tilde {\mathfrak R}  
\end{pmatrix},\quad \hat P_+= \begin{pmatrix}
   a^*+b^*\tilde {\mathfrak R}& b\eul^{-2\ii k\theta}\eul^{\frac{\ii  t}{2 k}}\\
0&a+b\tilde {\mathfrak R}^*  
\end{pmatrix}.
\end{align}
Then, using \eqref{calR} and \eqref{propR5}, straightforward calculations show \eqref{per_jump}.

\noindent (4) For $k\in  \{ \Im k > 0, k\in\Sigma_j\}$, we have
\begin{align}
\hat P_-^{-1}&=        
\begin{pmatrix}
  a+b\mathfrak R^*_-  &0\\
 b^*\eul^{2\ii k\theta}\eul^{-\frac{\ii  t}{2 k}}& a^*+b^*\mathfrak R_-
\end{pmatrix},\quad \hat P_+= \begin{pmatrix}
   a^*+b^*\mathfrak R_+&0\\
 -b^*\eul^{2\ii k\theta}\eul^{-\frac{\ii  t}{2 k}}&a+b\mathfrak R^*_+  
\end{pmatrix}.
\end{align}
Then, using \eqref{propR6} and \eqref{propR9} straightforward calculations show \eqref{per_jump}.

Similarly, one can show that \eqref{per_jump} holds in $\mathbb{C}^-$.
    
\end{proof}

\begin{proposition}\label{prop:beh_infty}
    As $k\to \infty$, $M^{(t)}(t,k)=I+O(\frac{1}{k})$ and $\hat M^{(t)}(t,k)=I+O(\frac{1}{k})$. 
\end{proposition}
\begin{proof}
     Notice that \eqref{s_inf} together with \eqref{spec_R} implies that
$(a+b\mathfrak R^*)(k)\to 1$ as $k\to\infty$ in $\mathbb{C}^+$, and \eqref{propR6} implies that $(a^*+b^*\mathfrak R)(k)\to1$ as $k\to\infty$ in $\mathbb{C}^+$. Furthermore, \eqref{s_inf} implies that $b^*\eul^{2\ii k\theta}\to 0$ as $k\to \infty$. Thus, the statement follows. 
\end{proof}

\begin{proposition}\label{prop:res}
(i) $M^{(t)(2)}(t,k)$ and $\hat M^{(t)(2)}(t,k)$ have simple poles at poles $\mu_j$ of $\mathfrak{R}$ with $\mu_j\in \mathbb{C}^+\setminus\{|k|<\epsilon\}$. Moreover,
\begin{align}\label{res_M_t_1}
\Res_{\mu_j} M^{(t)(2)}(t,k)& = M^{(t)(1)}(t,\mu_j)\ee^{\frac{\ii t}{2\mu_j}}\Res_{\mu_j}\mathfrak R, \quad \mu_j\in \mathbb{C}^+\setminus\{|k|<\epsilon\},
\end{align}
and $\hat M^{(t)(2)}(t,k)$ satisfies the same conditions 
\eqref{res_M_t_1}.

\noindent (ii)
$M^{(t)(1)}(t,k)$ and $\hat M^{(t)(1)}(t,k)$ have simple poles at poles $\bar\mu_j$ of $\mathfrak{R}^*$ with $\bar\mu_j\in \mathbb{C}^-\setminus\{|k|<\epsilon\}$. Moreover,
\begin{align}\label{res_M_t_3}
\Res_{\bar\mu_j} M^{(t)(1)}(t,k) = -M^{(t)(2)}(t,\bar\mu_j)\ee^{-\frac{\ii t}{2\bar\mu_j}}\Res_{\bar\mu_j}\mathfrak R^*, \quad \bar\mu_j\in\mathbb C^-\setminus\{|k|<\epsilon\},
\end{align}
and $\hat M^{(t)(1)}(t,k)$ satisfies the same conditions \eqref{res_M_t_3}.

\noindent (iii)
There are no  other pole singularities.    
\end{proposition}

\begin{proof}

If  $\mu_j$ is such that  $\mathfrak R(\mu_j)=\infty$, then $b^*(\mu_j)=0$. Dividing \eqref{calR} by $\mathfrak R$ and passing to the limit as $k\to \mu_j$, we obtain
\begin{equation}\label{Rb*}
  \lim_{k\to \mu_j} \mathfrak R(k)b^*(k)=-a^*(\mu_j)+a(\mu_j)\eul^{-2\ii \mu_j \theta}.  
\end{equation}
Consequently, 
\begin{equation}\label{a*-b*R}
(a^*+b^*\mathfrak R)(\mu_j)=a(\mu_j)\eul^{-2\ii\mu_j\theta}.  
\end{equation}
Now, combining \eqref{propR6}, \eqref{a*-b*R} and \eqref{det_rel}, we conclude that 

\begin{equation}\label{a-bR*}
(a+b\mathfrak R^*)(\mu_j)=a^*(\mu_j)\eul^{2\ii\mu_j\theta}.  
\end{equation}

Similarly, if $\bar \mu_j$ is such that $\mathfrak R^*(\bar \mu_j)=\infty$, then $b(\bar \mu_j)=0$. Proceeding as above,
we have
\begin{equation}\label{Rb*_2}
  \lim_{k\to \bar\mu_j} \mathfrak R^*(k)b(k)=-a(\bar\mu_j)+a^*(\bar\mu_j)\eul^{2\ii \bar\mu_j \theta}.  
\end{equation}
Consequently, 
\begin{equation}\label{a-bR*_2}
(a+b\mathfrak R^*)(\bar\mu_j)=a^*(\bar\mu_j)\eul^{2\ii\bar\mu_j\theta}  
\end{equation}
and
\begin{equation}\label{a*-b*R_2}
(a^*+b^*\mathfrak R)(\bar\mu_j)=a(\bar\mu_j)\eul^{-2\ii\bar\mu_j\theta}.  
\end{equation}

Recall that for $k\in \mathbb{C}^+\setminus\{|k|<\epsilon\}$, we have
\begin{align}\label{res_M_hat_M_1}
    & M^{(t)(1)}(t,k)=\hat M^{(1)}(0,t,k), \\\label{res_M_hat_M_2}
    & M^{(t)(2)}(t,k)=\hat M^{(2)}(0,t,k) ,\\\label{res_M_hat_M_3}
    & \hat M^{(t)(1)}(t,k)= (a^*+b^*\mathfrak R)(k)\hat M^{(1)}(\theta,t,k) -b^*(k)\ee^{2\ii k\theta}\ee^{-\frac{\ii  t}{2k}}\hat M^{(2)}(\theta,t,k),\\ \label{res_M_hat_M_4}
    & \hat M^{(t)(2)}(t,k)=(a+b\mathfrak R^*) \hat M^{(2)}(\theta,t,k).
\end{align}
First, let us consider the first columns. It is immediate that $M^{(t)(1)}(t,k)$ is regular at $\mu_j$. Moreover, combining \eqref{res_hatM-a} and \eqref{res_M_hat_M_3}, we obtain 
\begin{equation}\label{hatM_lim}
    \lim_{k\to \mu_j} \hat{M}^{(t)(1)}(t,k)=\hat{M}^{(1)}(\theta,t,\mu_j) a^*(\mu_j)
\end{equation}
and thus $\hat{M}^{(t)(1)}(t,k)$ is not singular at 
  $\mu_j$.

Considering the second columns, \eqref{res_hatM-a} and \eqref{res_M_hat_M_2} yield \eqref{res_M_t_1} for $M^{(t)}(t,k)$ while \eqref{res_hatM-a}, \eqref{a-bR*}, \eqref{hatM_lim} and \eqref{res_M_hat_M_4} give \eqref{res_M_t_1}
for $\hat M^{(t)}(t,k)$.

The singularities in $\mathbb{C}^-$  can be treated in a similar way.

\end{proof}

\begin{proposition}\label{prop:beh_i_2}
As $k\to0$ in either half-plane,
\begin{equation}\label{beh:P_0}
 \hat P(t,k)=I-\ii k(\theta-L)\sigma_3+O(k^2).
\end{equation}
In particular, $\hat\alpha(0,t)=\hat\alpha(\theta,t)$, $\hat\beta(0,t)=\hat\beta(\theta,t)$, and $\hat f_2(0,t)=\hat f_2(\theta,t)$.
\end{proposition}

\begin{proof}
The periodicity of $u_0$ and the integral equation at $t=0$ give
$\tilde a(k)=1+O(k^2)$ and $\tilde b(k)=O(k^2)$. Combining this with \eqref{s_til_s}, we have
\[
 a(k)=1+\ii k(\theta-L)+O(k^2),\qquad b(k)=O(k^2).
\]
The branch condition at the origin gives $\tilde{\mathfrak R}(k)=O(k)$: in \eqref{caltilR}, the coefficient of $\tilde{\mathfrak R}$ is $2\ii Lk+O(k^2)$ and $\tilde b=O(k^2)$.
Thus the diagonal entries of $\hat P$ are $1-\ii k(\theta-L)+O(k^2)$ and $1+\ii k(\theta-L)+O(k^2)$. The off-diagonal entry in $\mathbb C^+$ is $O(k^2)\eul^{-\ii t/(2k)}$, and that in $\mathbb C^-$ is $O(k^2)\eul^{\ii t/(2k)}$. These exponentials are bounded in their respective half-planes for $t\ge0$, proving \eqref{beh:P_0}.

The jump, normalization, residue, and boundedness conditions show by uniqueness that
$\hat M(0,t,k)=\hat M(\theta,t,k)\hat P(t,k)$.
Comparing \eqref{M_at_0_exp} and \eqref{beh:P_0} gives equality of $\hat\alpha$, $\hat\beta$, and $\hat f_2$ at the two endpoints. The diagonal linear term changes only $\hat f_1$.
\end{proof}

Consequently,
\begin{equation}
    \hat u(0,t)=\hat u(\theta,t), \quad  \hat w(0,t)=\hat w(\theta,t), \quad  \hat q(0,t)=\hat q(\theta,t).
\end{equation}
For $0\le y\le\theta$, the same RH comparison applies to $\hat M(y,t,k)$ and $\hat M(y+\theta,t,k)\hat P_y(t,k)$, where
\[
 \hat P_y(t,k)=\eul^{-\ii ky\sigma_3}\hat P(t,k)\eul^{\ii ky\sigma_3}.
\]
The jump and residue calculations are the same after conjugation by the spatial phase. The off-diagonal large-$k$ factors remain bounded for $y\ge0$, and the small-$k$ comparison is unchanged to first order off the diagonal. Thus the reconstructed coefficients agree at $y$ and $y+\theta$, and their periodic extension is $\theta$-periodic wherever the regular reconstruction exists.

We can now complete the verification of the initial condition. The equality $\hat w(y,0)=\hat w_0(y)$ and the normalization $x(0,0)=0$ give $u(x,0)=u_0(x)+C$. By Remark \ref{rem:mean-zero-RH}, the reconstruction has zero mean. Assumption \eqref{mean-zero-data} therefore gives $C=0$. Thus $u(x,0)=u_0(x)$.

\begin{theorem}\label{thm:periodic-reconstruction}
Let $u_0$ be real-valued, sufficiently smooth, and $L$-periodic.
Assume \eqref{mean-zero-data} and $u_{0x}(0)=u_{0x}(L)=0$.
Suppose its spectral data satisfy the nondegeneracy conditions stated
immediately before Proposition~\ref{prop:R}, and suppose the master
RH solution satisfies the reconstruction regularity conditions
stated after \eqref{res_hatM__loc}. On any time interval on which $\hat\alpha^2-\hat\beta^2>0$, the formulas \eqref{hmbbgg} and \eqref{x} reconstruct a smooth solution of \eqref{SP_1}, periodic with period $L$, satisfying $u(x,0)=u_0(x)$. In particular, if the indicated RH coefficients depend continuously on time at $t=0$, positivity at $t=0$ gives such an interval locally.
\end{theorem}

\begin{proof}
The reconstruction coefficient identities give the full system \eqref{SP_in_y}, including \eqref{SP_in_y-4}. The positive quantity $x_y=1/\hat q=\hat\alpha^2-\hat\beta^2$ makes the reciprocal map locally invertible, and \eqref{x} gives $x_t=-\hat u^2/2$. Proposition \ref{prop:hatu_u} then gives \eqref{SP}, hence \eqref{SP_1}. The periodicity argument and Proposition \ref{prop:periodic2} give period $L$. The preceding argument gives the prescribed initial value. At $t=0$, $1/\hat q_0>0$ has a positive minimum on the compact period, so continuity preserves positivity for sufficiently small time.
\end{proof}

\section{Acknowledgment}
This research was partially funded by the Research Council of Norway under project 361083 (Ukraina-NASTRAN cooperation).

\bibliographystyle{RS}
\bibliography{shepelsky_etal}

@article{IST01,
  title = {{Method for solving the Korteweg-de Vries equation}},
  author = {Gardner, Clifford S. and Greene, John M. and Kruskal, Martin D. and Miura, Robert M.},
  journal = {Physical Review Letters},
  volume = {19},
  issue = {19},
  pages = {1095--1097},
  numpages = {0},
  year = {1967},
  month = {Nov},
  publisher = {American Physical Society},
}

@article{IST02,
  title={{Exact theory of two-dimensional self-focusing and one-dimensional self-modulation of waves in nonlinear media}},
  author={V. E. Zakharov and A. B. Shabat},
  journal={Journal of Experimental and Theoretical Physics},
  year={1972},
  volume={34},
  pages={62--69}
}

@book{novikov1984theory,
  title={{Theory of solitons: the inverse scattering method}},
  author={Novikov, S. and Manakov, S. V. and Pitaevskii, L. P. and Zakharov, V. E.},
  year={1984},
  publisher={Springer Science \& Business Media}
}

@book{ablowitz1991solitons,
  title={Solitons, nonlinear evolution equations and inverse scattering},
  author={Ablowitz, M. J. and Clarkson, P. A.},
  volume={149},
  year={1991},
  publisher={Cambridge University Press}
}

@article{ablowitz1974inverse,
  title={The inverse scattering transform-{F}ourier analysis for nonlinear problems},
  author={Ablowitz, M. J. and Kaup, D. J. and Newell, A. C. and Segur, H.},
  journal={Studies in Applied Mathematics},
  volume={53},
  number={4},
  pages={249--315},
  year={1974},
  publisher={Wiley Online Library}
}

@book{trogdon2015riemann,
author = {Trogdon, Thomas and Olver, Sheehan},
title = {{Riemann--Hilbert problems, their numerical solution, and the computation of nonlinear special functions}},
publisher = {Society for Industrial and Applied Mathematics},
year = {2015},
address = {Philadelphia, PA}
}

@book{fokas2008unified,
author = {{F}okas, A. S.},
title = {A unified approach to boundary value problems},
publisher = {Society for Industrial and Applied Mathematics},
year = {2008},
}

@article{DFL21,
  title={The implementation of the unified transform to the nonlinear {S}chr{\"o}dinger equation with periodic initial conditions},
  author={Deconinck, B. and Fokas, A. S. and Lenells, J.},
  journal={Letters in Mathematical Physics},
  volume={111},
  number={1},
  pages={1--18},
  year={2021},
  publisher={Springer}
}

@article{FL21,
  title={A new approach to integrable evolution equations on the circle},
  author={Fokas, A. S. and Lenells, J.},
  journal={Proceedings of the Royal Society A},
  volume={477},
  number={2245},
  pages={20200605},
  year={2021},
  publisher={The Royal Society Publishing}
}

@article{SKPB24,
  title={{Periodic finite-band solutions 
to the focusing nonlinear Schr\"{o}dinger equation by the Fokas method: inverse and direct problems}},
  author={Shepelsky, D. and Karpenko, I. and Bogdanov, S.
          and Prilepsky, J.},
  journal={Proceedings of the Royal Society A},
  volume={480},
  pages={20230828},
  year={2024},
  publisher={The Royal Society Publishing}
}

@article{fokas1997,
  title={A unified transform method for solving linear and certain nonlinear {PDEs}},
  author={Fokas, A. S.},
journal={Proceedings of the Royal Society of London. Series A: Mathematical, Physical and Engineering Sciences},
  volume={453},
  pages={1411--1443},
  year={1997}
}

@ARTICLE{Z89,
  author={Zhou, X.},
  journal={SIAM Journal on Mathematical Analysis}, 
  title={{The Riemann–Hilbert Problem and Inverse Scattering}}, 
  year={1989},
  volume={20},
  pages={966-986}
}

@article{BSL17,
 author = {A.~{Boutet de Monvel} and D. Shepelsky  and L.  Zielinski},
 title = {{The short pulse equation by a Riemann–Hilbert approach}},
 journal = {Letters in Mathematical Physics},
 year = {2017},
 volume = {107},
 pages = {1345--1373},
}

@ARTICLE{SW04,
  author={Schäfer, T. and Wayne, C. E.},
  title={{Propagation of ultra-short optical pulses in nonlinear media}},
  journal={Physica D},
  volume={196},
  year={2004},
  pages={90--105}
}

@ARTICLE{CJWS05,
  author={Chung, Y. and Jones, C. K. R. T. and Schäfer, T. and Wayne, C. E.},
  title={{Ultra-short pulses in linear and nonlinear media}},
  journal={Nonlinearity},
  volume={18},
  number={3},
  year={2005},
  pages={1351--1374}
}

@ARTICLE{SS05,
  author={Sakovich, A. and Sakovich, S.},
  title={{The short pulse equation is integrable}},
  journal={J. Phys. Soc. Japan},
  volume={74},
  year={2005},
  pages={239--241}
}

@ARTICLE{BKSZ25,
  author={{Boutet de Monvel}, A. and Karpenko, I. and Shepelsky, D. and L.  Zielinski},
  title={{The periodic Camassa-Holm equation by the Riemann-Hilbert problem approach}},
  journal={arXiv:2511.20442v1},
year={2025}
}

@ARTICLE{KS25_,
  author={Karpenko, I. and Shepelsky, D.},
  title={{The Short Pulse equation on the half lines: a Riemann--Hilbert approach}},
  journal={arXiv:2512.05274},
year={2025}
}

@ARTICLE{Gambino2015,
  author={Gambino, G. and Tanriver, U. and Guha, P. and Choudhury, A. and Choudhury, S.},
  title={{Regular and singular pulse and front solutions and possible isochronous behavior in the short-pulse equation: phase-plane, multi-infinite series and variational approaches}},
  journal={Commun. Nonlinear Sci. Numer. Simul.},
  volume={20},
  number={2},
  year={2015},
  pages={375--388}
}

@ARTICLE{Matsuno2007,
  author={Matsuno, Y.},
  title={{Multiloop soliton and multibreather solutions of the short pulse model equation}},
  journal={J. Phys. Soc. Japan},
  volume={76},
  number={8},
  year={2007},
  pages={084003}
}

@ARTICLE{Sakovich2006,
  author={Sakovich, A. and Sakovich, S.},
  title={{Solitary wave solutions of the short pulse equation}},
  journal={J. Phys. A: Math. Gen.},
  volume={39},
  year={2006},
  pages={L361--L367}
}

@ARTICLE{Matsuno2008,
  author={Matsuno, Y.},
  title={{Periodic solutions of the short pulse model equation}},
  journal={J. Math. Phys.},
  volume={49},
  number={7},
  year={2008},
  pages={073508}
}

@ARTICLE{Coclite2015,
  author={Coclite, G. M. and di Ruvo, L.},
  title={{Well-posedness results for the short pulse equation}},
  journal={Z. Angew. Math. Phys.},
  volume={66},
  year={2015},
  pages={1529--1557}
}

@ARTICLE{Pelinovsky2010,
  author={Pelinovsky, D. and Sakovich, A.},
  title={{Global well-posedness of the short-pulse and sine-Gordon equations in energy space}},
  journal={Comm. Partial Differential Equations},
  volume={35},
  number={4},
  year={2010},
  pages={613--629}
}

@ARTICLE{Liu2009,
  author={Liu, Y. and Pelinovsky, D. and Sakovich, A.},
  title={{Wave breaking in the short-pulse equation}},
  journal={Dyn. PDE},
  volume={6},
  year={2009},
  pages={291--310}
}

@ARTICLE{F96,
  author  = {Fuchssteiner, B.},
  title   = {{Some tricks from the symmetry-toolbox for nonlinear equations: generalizations of the Camassa--Holm equation}},
  journal = {Phys. D},
  volume  = {95},
  year    = {1996},
  pages   = {229--243}
}

@ARTICLE{FO81,
  author  = {Fuchssteiner, B. and Fokas, A. S.},
  title   = {{Symplectic structures, their B{\"a}cklund transformations and hereditary symmetries}},
  journal = {Phys. D},
  volume  = {4},
  year    = {1981},
  pages   = {47--66}
}

@ARTICLE{OR96,
  author  = {Olver, P. J. and Rosenau, P.},
  title   = {{Tri-Hamiltonian duality between solitons and solitary-wave solutions having compact support}},
  journal = {Phys. Rev. E},
  volume  = {53},
  year    = {1996},
  pages   = {1900--1906}
}

@ARTICLE{Q06,
  author  = {Qiao, Z.},
  title   = {{A new integrable equation with cuspons and W/M-shape-peaks solitons}},
  journal = {J. Math. Phys.},
  volume  = {47},
  year    = {2006},
  pages   = {112701}
}

\end{document}